\documentclass[smallextended,envcountsect]{svjour3} 
\smartqed 

\usepackage{graphicx}
\usepackage{amsmath}
	\usepackage{lineno, blindtext}
\usepackage{enumitem}
\usepackage{amsfonts}

\usepackage{color}
\DeclareMathOperator{\grafo}{gph}
\DeclareMathOperator{\inte}{int}
\DeclareMathOperator{\supp}{supp}

\DeclareMathOperator{\co}{conv}
\DeclareMathOperator{\cco}{\cl{conv}}

\DeclareMathOperator{\epi}{epi}
\DeclareMathOperator{\dom}{dom}
\DeclareMathOperator{\cl}{cl}
\DeclareMathOperator{\sub}{\partial}

\newcommand{\R}{\mathbb{R}}
\newcommand{\N}{\mathbb{N}}
\newcommand{\Rex}{\overline{\mathbb{R}}}

\journalname{}
\let\epsilon\varepsilon

\begin{document}
\title{Formulae for the  Conjugate and the Subdifferential  of the  Supremum Function}

\subtitle{}

\author{Pedro P\'erez-Aros}

\institute{Pedro P\'erez-Aros \at
             Instituto de Ciencias de la Ingenier\'ia - Universidad de O'higgins\\
              Rancagua, Chile\\
              pedro.perez@uoh.cl \\
              This work was  supported in part by  CONICYT-PCHA/doctorado Nacional/2014-21140621.       
}

\date{Received: date / Accepted: date}

\maketitle

\begin{abstract}
The aim of this work is to provide formulae for the subdifferential and the conjungate function of the supremun function over an arbitrary family of functions. The work is principally motivated by  the case when  data functions are  lower semicontinuous proper and convex. Nevertheless, we explore the case when the family of functions is arbitrary, but satisfying that the biconjugate of the supremum functions is equal to the supremum of the biconjugate  of the data functions. The study focuses its attention on functions defined in finite-dimensional spaces, in this case the formulae can be simplified under certain qualification conditions. However, we show how to extend these results to arbitrary locally convex spaces without any qualification condition.
\end{abstract}


\section{Introduction}

	Convex analysis has been one of the most studied topics in nonsmooth analysis; currently, chapters  referring to convex analysis can be found in many monographic books related to optimization and variational  analysis. The development of convex analysis has brought many tools to establish necessary and sufficient conditions for optimality. A remarkable tool in this scenario is the notion of the subdifferential of a convex function. Everybody interested in nonsmooth analysis is aware of classical calculus rules and qualification conditions, together with practical algorithms, which involve this mathematical object.
	
		Another important tool in (convex) subdifferential  calculus is the so called  {\emph{approximate subdifferential of convex functions}} (see, e.g., \cite{MR1261420}). This operator has played a central role in optimization; its impact can be seen in applications such as  linear programming, convex and non-convex programming, stochastic programming and semi-infinite programming among many other topics, because this mathematical object has been used to provide calculus rules without qualification conditions (see, e.g., \cite{MR1245600,MR1330645}),  characterize points of  {sub-optimality} (see, e.g., \cite{MR1921556}). It has been used in many minimization algorithms (see, e.g., \cite{MR1261420}) and also can be found in connection with integration of subdifferentials (see e.g \cite{MR2155565,MR2861329,CORREA2018}). Also, the nonemptiness of this operator is guaranteed in arbitrary locally convex spaces under weaker conditions (\cite{MR1921556}) and it represents, in some sense, a continuous multifunction (see, e.g., \cite{MR600082}). All of these properties have motivated some authors to give generalizations of this  {object to vector optimization}  (see, e.g., \cite{MR3243824,MR3005026} and the references therein).
		
		A functional, which appears commonly in applications, is the supremum functional of a family of functions, for this reason many authors have shown several formulae for the calculus of the  subdifferential of the supremun function (see, e.g.,  \cite{MR2489616,MR3561780,MR2448918,MR2837551,MR3033113,2018arXiv180100724P}  and the reference therein). Moreover in \cite{MR2448918} the   authors have shown that the supremum function appears to be fundamental in the theory of convex subdifferential calculus,  in the sense that the general formula presented in \cite[Theorem 4]{MR2448918} allows the authors to recover classical and important formulae available in the literature (see \cite[Corollary 16]{MR2448918}, where the authors recover as a simple corollary of their general formula the calculus rules for the sum and composition). However, none of them have studied the $\varepsilon$-subdifferential of the supremum function. This observation motivates this work to provide calculus for the subdifferential  together with formulae for the conjugate function of the supremum function. 
		
		This work is organized as follows. In Section \ref{Notation}, we present some classical  notation and definitions of convex analysis which agree with many monograph books (see, e.g., \cite{MR1921556,MR0467080,MR1451876}).  In Section \ref{Prelimiliminarresult}, we give preliminary results concerning the calculus of the conjugate and the subdifferential of the  supremum function.  Section \ref{ORDEREDSETS} is devoted to studying formulae for the conjugate function of the supremun function of an arbitrary  family of (possibly non-convex) functions using the additional assumption that the index set is an ordered set and the functions  are epi-pointed; this allows us to give formulae for the conjugate of the $\limsup$ and $\liminf$ function. Finally, we divide Section \ref{SECTION:Calculus}, where the main results are  {established}, in two subsections. Subsection \ref{FINITEDIMENSION} is focused on when the  functions are defined in a  finite-dimensional space; in this scenario we investigate the subdifferential  and the conjugate function of the pointwise supremum under the standard assumption that the normal cone of the domain of the supremum function at the point of interest  does not contain lines.  The main results about the calculus of the subdifferential of the supremum function is Theorem  \ref{TEO5.2}; later under assumptions of compactness of the index set and some continuity property, which are classical in semi-infinite programming, we derive Theorem \ref{teorem:ref}. Posteriorly  we get a result that can be understood as the conjugate counterpart of the last two theorems (see Theorem \ref{THEOREMCONJUGATE}).  Due to the fact that this research is principally motivated by \cite{MR2448918}, in this section  we bypass the convexity of the data function  by the weaker assumption that the biconjugate of the  supremum function is equal to the supremum of the biconjugate of the data functions. This kind of hypothesis has been recently used in several works (see, e.g., \cite{MR3561780,MR2448918,MR2837551,MR2731292,MR2783793}). Subsection \ref{INFINITEDIMENSION} is motivated by the completeness of this work {. In} this subsection we show how to generalize the  mentioned results (given in a  finite-dimensional space) to a general locally convex space and without any qualification condition using the family of all  finite-dimensional subspaces containing a given point. 

\section{Notation}\label{Notation}
Throughout the paper $X$ and $X^\ast$ will denote  two (separated) locally convex spaces (lcs) in duality by the bilinear form $\langle \cdot,\cdot \rangle : X^\ast\times X  \to \mathbb{R}$, given by  $\langle x^\ast, x \rangle =x^\ast(x)$. In  $X^\ast$ the weak topology is denoted by  $w(X^\ast,X)$  ($w^\ast$ for short) and  the Mackey topology  on $X^\ast$ is  denoted by  $\tau(X^\ast,X)$, the space $X$ will be endowed with a compatible initial topology $\tau$ (i.e., $(X,\tau)^\ast=X^\ast$). The family of all closed convex and balanced neighborhoods of zero for the topology $\tau$ will be denoted by $\mathcal{N}_0(\tau)$ and we omit the symbol $\tau$ when there is no confusion, similar terminology is adapted for the space $X^\ast$.   We will write $\overline{\mathbb{R}}:= [-\infty,\infty ]$ and we adapt the convention $\frac{\alpha}{0}=+\infty$ for all $\alpha>0$. The symbol $\dim X$ means the dimension of the space $X$.

{ {The closure of $A$ will be denoted by $\cl A$. We denote by $\inte(A)$,  $\textnormal{conv}(A)$ and  $\cco(A)$, the interior, the \emph{convex hull} and the \emph{closed convex  hull} of $A$, respectively}}. The \emph{polar} of $A$ is the set  $A^o := \{ x^\ast\in X^\ast : \langle x^\ast,x\rangle \leq 1, \forall x\in A\}$.
When the space $X$ is  finite-dimensional, the norm on $X$ and $X^\ast$ will be denoted by $\| \cdot \|$ and $\| \cdot \|_\ast$, respectively. Given $x \in X$ (or $x^\ast\in X$) and $r\geq 0$ we denote $\mathbb{B}(x,r):=\{y \in X: \|x-z\|\leq r  \}$ ($\mathbb{B}(x^\ast,r):=\{y^\ast \in X^\ast: \|x^\ast-z^\ast\|_{*}\leq r  \}$) .

For a given function $f:X\to\overline{\mathbb{R}}$, the (effective) \emph{domain} and the  \emph{epigraph}  of $f$ are $$\dom f := \{ x\in X: f(x) < +\infty \}  \text{ and } \epi f := \{ (x,\lambda)\in X\times \R : f(x)\le\lambda \},$$ respectively. We say that $f$ is \emph{proper} if $\dom f \neq\emptyset$ and $f>-\infty$, and  \emph{inf-compact} if for every $\lambda \in \mathbb{R}$ the set $[f\le \lambda]:=\{ x \in X: f(x)\leq \lambda \}$ is compact. We denote $\Gamma_0(X)$  the class of proper lower semicontinuous (lsc) convex functions on $X$. The \emph{Fenchel conjugate} of $f$ is the function $f^\ast:X^\ast \to\overline{\mathbb{R}}$ defined by
\[ f^\ast(x^\ast) := \sup_{x\in X}\{ \langle x^\ast,x\rangle - f(x)\}, \]
and the \emph{biconjugate}  of $f$ is $f^{\ast \ast}:= (f^\ast)^\ast:X \to \Rex$.  {The \emph{closed hull}, the \emph{convex hull} and   the \emph{closed convex hull} of  $f$ are denoted  by  $\cl{f}$, $\co{f}$ and  $\cco f$ respectively, and they are defined  as the functions such that $\epi \cl{f}=\cl \epi f $, $\epi \co{f}=\co \epi f $,  and  $\epi \cco f = \cco\epi f$ respectively. We recall that $f^{\ast \ast}\leq \cco f$ and whenever $f^{\ast \ast}$ is proper one has $f^{\ast \ast}= \cco f$.} 

For $\varepsilon\geq0$, the \emph{$\varepsilon$-subdifferential} (or the \emph{approximate subdifferential of convex functions})   of $f$ at a point $x\in X$, where it is finite, is  the set
\[ \partial_\varepsilon f(x) := \{ x^\ast\in X^\ast : \langle x^\ast,y-x\rangle \leq f(y) - f(x) +\varepsilon,\ \forall y\in X  \}; \]
if $f(x)$ is not finite, we set $\partial_\varepsilon f(x):=\emptyset$. The special case $\varepsilon=0$ is the classical convex subdifferential, also called  the  \emph{Moreau-Rockafellar Subdifferential}  and  is denoted by $\partial f(x)$. 

The \emph{indicator} and the \emph{support} functions of a set $A$ ($\subset X,X^\ast$) are, respectively,
\[ \delta_A(x) := \begin{cases}
0,\qquad & x\in A,\\
+\infty,& x\notin A,
\end{cases}\qquad\qquad \sigma_A :=\delta_A^\ast. \]
The \emph{asymptotic cone} of $A$ is defined by $A_\infty:=\bigcap_{\varepsilon >0 } \cl \big(   ] 0,\varepsilon] A \big)$ (see, e.g., \cite{MR562254}). When the space $X$ is a Banach space the   asymptotic cone is commonly  defined as $A_\infty= \{ d \in X : \exists s_k \searrow 0 \text{ and } x_k \in A \text{ such that } \lim \frac{x_k}{s_k} = d\}$. When the set  $A$ is convex $A_\infty$ is also known as the \emph{recession cone}; in this case it admits the following representation $$A_\infty= \{ d \in X:   {x+ \lambda d} \in A \text{ for some }x \text{ and all }\lambda \geq 0  \}.$$ We set 
\begin{align*}
A^{-}&:=\{ x^\ast\in  X^\ast: \langle x^\ast , x \rangle \leq 0,\; \forall x\in A   \},\\
A^{\perp}&:=\{ x^\ast\in  X^\ast: \langle x^\ast , x \rangle = 0,\; \forall x\in A   \}
\end{align*} the \emph{negative dual cone}, and the \emph{orthogonal subspace} (or
annihilator) of $A$, respectively.

The $\varepsilon$-normal set of $A$ at a point $x$ is defined and denoted by $$ {\textnormal{N}^\varepsilon_A(x):=\partial_{\varepsilon} \delta_A(x).}$$

Now let us introduce the class of epi-pointed functions. This class of functions has shown that it shares  important properties with the class of convex function, we refer to \cite{MR2861329,CORREA2018,MR3507100,Correa2017,MR1884906,MR3509670,INTECONV,INTECONV2,2018arXiv180107378P}  and the references therein.

	\begin{definition}\label{definitionepipointed}
	A function $f:X\to \Rex$ is said to be epi-pointed if $f^\ast$ is proper and  $\tau(X^\ast,X)$-continuous at some point of its domain. 
\end{definition}

 {As far as we know, this family of functions was introduced in finite dimensions in \cite{MR1416513}. However, the the above extension was introduced in \cite{MR3509670} with the name of \emph{Mackey epi-pointed functions}.}

In a finite-dimensional Banach space $X$ the reader can understand the above property in terms of the pointedness of the epigraph of  its \emph{asymptotic function}, which is  given by 
$$ f^\infty (u):=\liminf\limits_{ \substack{ s \searrow 0\\ w \to u }} s f( s^{-1} w	) 	, $$
or in terms of a  \emph{coercive property}, which is the  existence of   $x^\ast \in X$, $\alpha >0$ and $r \in \R$ such that $f(x) \geq \langle x^\ast , x \rangle +\alpha\| x\| +r $ for all $x\in X$  (see, e.g., \cite[Proposition 3.1]{Correa2017}). 

For a set $T$, $\mathcal{P}_f(T)$  denotes the set of all $F\subseteq T$ such that $F$ is finite. We  define $$f(x):=\sup_{t \in T} f_t,\; \text{ for }x\in X$$ as the pointwise supremum of the family $\{f_t : t \in T \}$. 	When $T$ is a directed set  ordered  by $\preceq$ (i.e.,  $(T,\preceq)$ is an ordered set and for every $t_1,t_2\in T$ there exists $t_3 \in T$ with $t_1\preceq t_3$ and $t_2 \preceq t_3$) we say that the family of functions is increasing provided that  for all $t_1,t_2 \in T$
\begin{align*}
t_1 \preceq t_2 \implies f_{t_1} \leq f_{t_2} \;\text{(i.e., }  f_{t_1}(x) \leq f_{t_2}(x) , \; \forall x\in X\text{)}.
\end{align*}   	
  For  convenience we also use the notation $t_2 \succeq t_2$ iff $t_1 \preceq t_2$.

Following the standard notation (see, e.g., \cite{MR2489616,MR3561780,MR2448918}), for a subset $Z$ of $\R$  let $Z^{(T)}$ be the set of all $(\lambda_t) \in Z^{T}$ such that $\lambda_t\neq0$ for finitely many $t\in  T$. The support of $\lambda$ is defined as $\supp \lambda :=\{t \in T : \lambda_t \neq 0  \}$. The \emph{generalized simplex} on $T$ is the set $\Delta(T):=\{\lambda \in  [0,1]^{(T)}  : \sum_{t\in T} \lambda_t =1 \}$, also it will be convenient to denote $ \Delta^\varepsilon(T): =\{\lambda \in  [0,\epsilon]^{(T)}  :  \sum_{t\in T} \lambda_t =\varepsilon\}$.  For a  family of functions $\{ f_t \}_{t\in T} \subseteq {\Rex}^T$,  a point $\bar{x}$ and $\varepsilon\geq 0$,  the set of $\varepsilon$-active indexes at $\bar{x}$ is defined by $T_\varepsilon(\bar{x}):=\{ t\in T : f_t(\bar{x}) \geq f(\bar{x}) -\epsilon		\}$.

\section{Preliminary Result}\label{Prelimiliminarresult}

\subsection{Basic Properties of the Conjugate and the $\varepsilon$-Subdifferential of the Pointwise Supremum}

The next lemma shows some basic  properties of the the pointwise supremum  function.
\begin{lemma}\label{Lemma2}
	Let $\{f_t : t \in T \}$ be an arbitrary family of functions and define $h:=\inf_{t \in T} f^\ast_t $. Then
	\begin{enumerate}[label={(\alph*)},ref={(\alph*)}]
		\item\label{Lemma2b}   $h^\ast=\sup f^{\ast \ast}_t$, consequently if $h^{\ast \ast}\neq -\infty$, then $\cco h=h^{\ast \ast}$. Therefore,  {if the functions $f_t \in \Gamma_0(X)$} and $f$ is proper one has $$\epi f^\ast= \cco \bigcup_{t\in T}  \epi f^\ast_t .$$
	
		\item \label{Lemma2CONVEXENVELOPE} $\co  h (x^\ast)=\co \{  f_t^\ast : t\in T  \}(x^\ast)$, where $$\co \{  f_t^\ast : t\in T  \}(x^\ast):= {\inf\Bigg\{  \sum_{t\in T} \lambda_t f_t(x^\ast_t) : \begin{array}{c}
	\lambda \in \Delta(T) \text{ and }\\ \sum_{t\in t}\lambda_t x^\ast_t =x^\ast 
		\end{array}    \Bigg\}}.$$ Moreover,  the infimum can be taken only in $\lambda \in \Delta(T)$ with $$\# \supp \lambda \leq  \min\{  \dim X + 1, \# T     \}. $$
			\item\label{Lemma2important2}   If  $f^\ast$ is proper and   $f^{\ast \ast}=\sup_{t\in T} f^{\ast \ast}_t$, then  $f^\ast=\cco\{ f_t^\ast: t\in T\}$.
		\item\label{Lemma2d}  If $\{f_t : t \in T \}$  is an increasing family of functions, then $\inf_{t\in T} f^\ast_t$ is convex.
		
		\item\label{Lemma2c} If $f_t \leq f_s$ and $f_t$ is epi-pointed and $f^\ast_s$ is proper, then $f_s$ is epi-pointed.
		\item\label{Lemma2important}  {For every $t\in T$ we have that  $\dom f_t^\ast \subseteq \dom f^\ast $ and $\inte \dom f_t^\ast \subseteq \inte \dom f^\ast $}. In addition, if $f^{*}$ is proper,  {$f^{\ast \ast}=\sup_{t\in T} f^{\ast \ast}_t$} and $\{f_t : t \in T \}$  is an increasing family of \mbox{epi-pointed} functions one has 
		\begin{align*}
		\bigcup_{ t\in T}\inte \dom f_t^\ast = \inte \dom f^\ast \text{ and } \bigcup_{ t\in T}\inte \epi  f_t^\ast = \inte \epi f^\ast .
		\end{align*}
	\end{enumerate}
\end{lemma}

{\it Proof}
	\begin{enumerate}[label={(\alph*)},ref={(\alph*)}]
		\item  {It is not difficult to see that  $h^\ast=\sup f^{\ast \ast}_t$, then $\cco h=h^{\ast \ast}$ whenever $h^{\ast \ast}$ is proper (see, e.g., \cite[Theorem 2.3.4]{MR1921556}).   Moreover, when $\{f_t \}_{t\in T} \subseteq \Gamma_0(x)$ we have $f_t^{\ast \ast} =f_t$ for each $t\in T$, in particular,  $h^\ast=f$. Then, using the fact that  $f^\ast$ is proper we have $h^{\ast \ast}$ is also proper and consequently  $$\epi f^\ast = \cco\bigcup_{t\in T}  \epi f^\ast_t \text{ and } f^\ast=\cco \inf_{t\in T} f_t^\ast.$$} 
		
		\item  { First, $\co h \leq \co\{  	 f^\ast_t : t\in T\}$ due to the convexity of $\co h$. Now, we notice that \begin{align}\label{stricepigraph}
			\epi_s h = \bigcup_{ t\in T} \epi_s f^\ast_t,
			\end{align}
			where $\epi_s$ denotes the strict epigraph. It follows that for every $x^\ast \in \dom \co h$ and every and   $\epsilon >0$, the element  $(x^\ast,\alpha_\epsilon) \in \co (\epi h)$, where 
			$$
		\alpha_\epsilon:=\Bigg\{ \begin{array}{cl}
			\co h(x^\ast), &\text{ if } h(x^\ast)\in\R, \\
	\epsilon^{-1}, &\text{ if } h(x^\ast)=-\infty.
			\end{array}$$ 
			Then,   there exists $\lambda_i \geq 0$ and $(x_i^*, \beta_i ) \in \epi h$ with $i=1,...,p$ and $\sum_{ i=1  }^p \lambda_i =1$ such that $(x^\ast,\alpha_\epsilon) = \sum_{ i=1  }^p \lambda_i (x_i^*, \beta_i )$. Moreover, by \eqref{stricepigraph} there exists $t_i \in T$  such that $(x_i^*, \beta_i + \epsilon ) \in \epi_s f^*_{t_i}$. It implies that
		\begin{align*}
		\alpha_\epsilon + \epsilon  \geq  \sum_{ i=1  }^p \lambda_{i} f^*_{t_i} (x^*_{i}) \geq \co \{ f_t^\ast : t\in T \}.
		\end{align*}
		Since $\alpha_\epsilon + \epsilon$ converges to $\co h(x^\ast)$ as $\epsilon \searrow 0$ we conclude the equality.
	 }	
 Moreover, when the space is  finite-dimensional we can use Carathéodory's Theorem to consider only $\lambda \in \Delta(T)$ with $\#\textnormal{supp} \lambda \leq \dim X + 1$.
		
		\item  {Since $f^\ast$ is proper and   $h^{\ast \ast} \geq f^\ast$ we have $h^{\ast \ast } \neq -\infty$, then by \ref{Lemma2b} we get  $h^{\ast \ast} = \cco h$. Finally, let us prove that $\cco h \leq  f^\ast$, indeed}
		\begin{align*}
	 {	h^{\ast\ast}}&= \sup_{x\in X}\{ \langle  \cdot ,x \rangle - h^{*}(x) \}=\sup_{x\in X}\{ \langle  \cdot ,x \rangle - f^{\ast \ast}(x) \} = f^\ast.
		\end{align*}
		Therefore $f^\ast= h^{\ast \ast}= \cco h= \cco \{ f_t^\ast : t \in T\}$.
		\item Because the family of functions is increasing, we have  $\epi f_t^\ast \subseteq \epi f_s^\ast$ for $t\leq s$, therefore $\bigcup_{ t\in T} \epi f^\ast_t $ is convex and then 
		$\inf_{t\in T} f_t^\ast$ is convex.
		\item If $f_t \leq f_s$ one has $f^\ast_t \geq f_s^\ast$ and since $f^\ast_t$ is bounded in $\inte(\dom f_t^\ast)$, so is $f_s^\ast$, which implies the continuity of $f_s^\ast$.
		\item First, one  has $\bigcup_{ t\in T}\dom f_t^\ast \subseteq \dom f^\ast $ and $ \bigcup_{t\in T} \epi f^\ast_t\subseteq \epi f^\ast$, consequently $\bigcup_{ t\in T}\inte \dom f_t^\ast \subseteq \inte \dom f^\ast $ and $ \bigcup_{t\in T} \inte \epi f^\ast_t \subseteq \inte \epi f^\ast$. Now, by \ref{Lemma2important2}  and \ref{Lemma2d}    one has 
		\begin{align}\label{0001}
		\epi f^\ast = \cl\{  \bigcup_{t\in T} \epi f^\ast_t\}=  \cl\{  \bigcup_{t\in T} \inte \epi f^\ast_t\}
		\end{align}
	 and then $ \dom f^\ast  \subseteq \cl\{   \bigcup_{ t\in T}\inte \dom f_t^\ast   \}$. Fix   $(u^\ast_0,\beta_0) \in  \bigcup_{t\in T} \inte \epi f^\ast_t$ and take  $(x_0^\ast, \alpha_0)  \in \inte \epi f^\ast$; then there exists   $\gamma >0$ such that $(x,\alpha):=(x_0^\ast,\alpha_0) + \gamma ( (x_0^\ast,\alpha_0) - (u^\ast_0,\beta_0) ) \in  \epi f^\ast$, because  $ \bigcup_{t\in T} \inte \epi f^\ast_t$ is a convex open  dense subset in $ \epi f^\ast$ (see \eqref{0001}) one has 
		$$\frac{1}{1+\gamma}  (x,\alpha) + \frac{\gamma }{1+\gamma} (u^\ast_0,\beta_0)= (x_0^\ast,\alpha_0) \in  \bigcup_{t\in T} \inte \epi f^\ast_t.$$
		Using the same argument one also can conclude  $\bigcup\limits_{ t\in T}\inte \dom f_t^\ast= \inte \dom f^\ast $.
	\end{enumerate}
  \qed

\begin{remark}
	  Lemma \ref{Lemma2} \ref{Lemma2CONVEXENVELOPE} appears to be very simple. Nevertheless, it represents a key point in  the entire work and it allows us to use a weaker condition than the lsc and convexity of the data function (i.e. $ f_t^\ast \in  \Gamma_0(X)$ for all $t\in T$). In this work we use the hypothesis that $f^{\ast \ast}= \sup_{t\in T}  f^{\ast \ast}_t$,  which has been used recently by some authors (see  \cite{MR2731292,MR2783793,MR2960099}) and it is weaker than the hypothesis used in \cite{MR2448918,MR2837551}. More precisely in \cite{MR2448918} the authors use the hypothesis that $\{f_t\}_{t\in T}$ are proper convex functions satisfying $\cl f =\sup_{t\in T} \cl f_t$. Posteriorly in \cite{MR2837551} the authors use the  hypotheses that $\{f_t\}_{t\in T}$ are proper convex functions and $\cl f(x) =\sup_{t\in T} \cl f_t(x)$ for all $x\in  \cl (\dom f)$. Recently, in \cite{MR3561780} the authors presented an improvement of the last results, in this work the functions  $\{f_t\}_{t\in T}$ are assumed to be proper convex functions and $\cl f(x) =\sup \cl f_t(x)$ for all $x\in  \dom f$ (see \cite[Corollary 6]{MR3561780}). It is worth mentioning that to prove Lemma \ref{Lemma2} \ref{Lemma2important2}  one can  assume the weaker hypothesis that $f^\ast$ is proper and  $ f^{\ast \ast}(x)= \sup_{t\in T}  f^{\ast \ast}_t(x)$ for all $x \in D$, where $D$ is a graphically dense subset of $\dom  f^{\ast \ast}$ (which is satisfied under the assumption of \cite[Corollary 6]{MR3561780}). For the sake of simplicity, we keep in mind  these sophisticated assumptions for a  future work. 
\end{remark}

The next lemma is a slight extension of \cite[Lemma 4.5]{Correa2017} without the assumption of convexity.

\begin{lemma}\label{Upper:stimationeta}
	Let $g:X \to \Rex$ be an epi-pointed function.  {Consider} $\eta >0$ and $x \in X$ such that  $\eta +g^{\ast \ast}(x)>g(x)$. Then
	\begin{align}\label{Equation:Subdifferential}
	\partial_\eta g(x)= \cl^{w^\ast}\bigg\{\partial_\eta g(x)\cap \inte \dom g^\ast   \bigg\}.
	\end{align}
Consequently  for every $\eta \geq 0$
\begin{align}\label{Equation:Subdifferential2}
	\partial_\eta g(x)= \bigcap\limits_{ \gamma >   0}\cl^{w^\ast}\bigg\{\partial_{\eta+\gamma} g(x)\cap \inte \dom g^\ast   \bigg\}.
\end{align}
\end{lemma}

{\it Proof} 	First, if $|g^{\ast \ast}(x)|=+\infty$ one has the empty set in both sides of (\ref{Equation:Subdifferential}). Now consider $g^{\ast \ast}(x)\in \R$ and take $\gamma >0$  such that $\eta-\gamma +g^{\ast \ast}(x)>g(x)$. Because $g^{\ast \ast}$ is a proper convex and lsc function and $\gamma>0$, one has $\partial_\gamma  g^{\ast \ast}(x)$ is non-empty, and by  \cite[Lemma 4.5]{Correa2017}  we get  $\partial_\gamma  g^{\ast \ast}(x) \cap  \inte \dom g^\ast \neq \emptyset$.

 Then taking $u^\ast \in\partial_\gamma  g^{\ast \ast}(x) \cap  \inte \dom g^\ast$ one gets $u^\ast \in \partial_\eta g(x)$, which implies, $\partial_{\eta} g(x) \cap \inte \dom g^\ast  \neq \emptyset$.   {Since $\inte \dom g^\ast$ is open, convex and dense in the convex set $\dom g^\ast $,} we have that
	$$ \cl \bigg\{\partial_\eta g(x)\cap \inte \dom g^\ast    \bigg\}= \partial_\eta g(x) \cap \cl \inte \dom g^\ast = \partial_\eta g(x).$$ Finally,  (\ref{Equation:Subdifferential2}) follows noticing that if $\partial_{\eta }g(x) \neq \emptyset$ we have $\gamma + \eta + g^{\ast \ast}(x) > g(x)$ and  $\partial_{\eta }g(x) = \bigcap_{ \gamma >   0} \partial_{\eta + \gamma}  g(x)$.
\qed     

\begin{lemma}\textnormal{\cite[Lemma 1.1]{MR1330645}}\label{LEMMA:HIRIART}
	Let $h$ be an extended-real-valued convex function defined over $X^\ast$. Then, for all $r \in \R$,
	\begin{align*}
	\{  x^\ast \in X^\ast :  \cl{h}(x^\ast)\leq r      \}=\bigcap\limits_{ \gamma >0 } \cl \{ 	x^\ast \in X^\ast : h(x^\ast) < r +\gamma	\}.
	\end{align*}
	Moreover, if $r > \inf_{X^\ast} h$, then 
	\begin{align*}
	\{  x^\ast \in X^\ast :  \cl{h}(x^\ast)\leq r      \}=\cl \{ 	x^\ast \in X^\ast : h(x^\ast) < r 	\}.
	\end{align*}
\end{lemma}
 {The following results give us a first representation of the $\epsilon$-subdifferential of $f$, this result corresponds to a slight extension  of  \cite[Theorem 2]{MR3740511}.}
\begin{proposition}\label{Epsilonformula}
	Let $\{f_t : t \in T \}$ be an arbitrary family of functions such that $f^{\ast \ast}=\sup_{T}  f^{\ast \ast}_t$. 	If  $\varepsilon > \inf_{x^\ast \in X^\ast} \{ f^\ast(x^\ast) +f(x) -\langle x^\ast ,x \rangle 			\} $, then
	\begin{align}\label{FORMULABASIS}
	\partial_\varepsilon f(x) =     \cl\bigg\{  \sum_{t\in \supp T} \lambda_t \partial_{\eta_t }  f_t(x) :  \begin{array}{c}
	\lambda\in \Delta(T), \; \eta_t \geq 0,\; \\ \sum_{t\in T}\lambda_t \cdot \eta_t \in [0,\varepsilon) \text{ and }\\
		\sum_{t\in T} \lambda_t \cdot  f_t  (x) > f(x) + \sum_{t\in T}\lambda_t \cdot \eta_t - \varepsilon 	
	\end{array}				\bigg\}.
	\end{align}
	Consequently, for every  $\varepsilon \geq 0$
	\begin{align}\label{FORMULABASIS2}
	\partial_{\varepsilon} f(x) = \bigcap\limits_{\gamma >   \varepsilon } \cl\bigg\{ \hspace*{-0.1cm} \sum\limits_{ t\in \supp T} \lambda_t \partial_{\eta_t }  f_t(x) : \hspace*{-0.1cm} \begin{array}{c}
	\lambda\in \Delta(T), \; \eta_t \geq 0,\; \\\sum_{t\in T}\lambda_t \cdot \eta_t \in [0,\gamma),
	\text{ and }\\	\sum\limits_{t\in T} \lambda_t  f_t  (x) \geq f(x) + \sum_{t\in T}\lambda_t \cdot \eta_t - \gamma 	
	\end{array}				\bigg\}.
	\end{align}
	
 Moreover, when the space $X$ is  finite-dimensional  {the infimum can be taken only} $\lambda \in \Delta(T)$ with  $\# \supp \lambda \leq  \min\{  \dim X + 1,\#T \}$.
\end{proposition}
{\it Proof}
	Since  the right side of (\ref{FORMULABASIS}) and (\ref{FORMULABASIS2}) are included in $\sub_\varepsilon f(x)$ we must focus on the opposite inclusion. Assume that 
	$\partial_\varepsilon f(x)\neq \emptyset$, then $f^\ast$ is proper and by Lemma \ref{Lemma2}\ref{Lemma2important2}, $f^\ast=\cco \{ f^\ast_t : t\in T \}$. Hence, using  Lemma \ref{LEMMA:HIRIART}
	$$\partial_{\varepsilon } f(x) = \cl\{   x^\ast \in X^\ast : \co \{ f_t^\ast : 	 t\in T\}(x^\ast) +f(x) < \langle  x^\ast , x \rangle + \varepsilon \}.$$
	Now let $x^\ast\in X^\ast$ satisfying $\co \{ f_t^\ast : 	 t\in  T\}(x^\ast) +f(x)  -\langle  x^\ast , x \rangle <   \gamma	$, then by Lemma \ref{Lemma2} \ref{Lemma2CONVEXENVELOPE} there exist $\lambda \in \Delta(T)$ (with $\#\textnormal{supp} \lambda \leq \dim X + 1$) and $x^\ast_t \in X^\ast$ for $t\in\supp \lambda$  such that $\sum_{t\in T} \lambda_t x^\ast_t  = x^\ast$ and 
	\begin{align*}
	\co \{ f_t^\ast : 	 t\in T\}(x^\ast) +f(x) - \langle  x^\ast , x \rangle  \leq&  \sum_{t\in T} \lambda_t \big(  f_t^\ast(x_t^\ast)  +f_t(x) - \langle x^\ast_t, x\rangle				 \big)  +\\ &\sum\lambda_t \big( f(x)- f_t(x) \big)		<  \varepsilon,
	\end{align*}
	then taking  $\eta_t := f_t^\ast(x_t^\ast)  +f_t(x) - \langle x^\ast_t, x\rangle $ if $t \in \supp \lambda$ and $\eta_t:=0$ if $t\notin  \supp \lambda$ one gets that $x^\ast_t \in \partial_{\eta_t } f_t(x)$ for every $t\in \supp \lambda$,  $\sum_{t\in T}\lambda_t \cdot \eta_t \in [0,\varepsilon)$ and $\sum_{t\in T} \lambda_t \cdot  f_t  (x) > f(x) + \sum_{t\in T}\lambda_t \cdot \eta_t - \varepsilon$, which concludes the proof of  (\ref{FORMULABASIS}).
	
	Finally, if  $\varepsilon \geq 0$ is any real number such that $\sub_\varepsilon f(x)\neq \emptyset$ one has that $\inf_{x^\ast \in X^\ast} \{ f^\ast(x^\ast) +f(x) -\langle x^\ast ,x \rangle 			\}\leq \varepsilon $, consequently using (\ref{FORMULABASIS}) for  $\gamma >\varepsilon$ and  the fact that $\partial_\varepsilon f(x) =   \bigcap_{\gamma >\varepsilon} \partial_\gamma f(x) $ we obtain
\begin{align*}
\partial_\varepsilon f(x) &= \bigcap_{\gamma >\varepsilon}   \cl\bigg\{  \sum_{t\in T} \lambda_t \partial_{\eta_t }  f_t(x) :  \begin{array}{c}
\lambda\in \Delta(T), \; \eta_t \geq 0,\;\\
 \sum_{t\in T}\lambda \cdot \eta_t \in [0,\gamma) \text{ and}\\
	\sum_{t\in T} \lambda_t \cdot  f_t  (x) > f(x) + \sum_{t\in T}\lambda_t \cdot \eta_t - \gamma 	
\end{array}				\bigg\}.
\end{align*}	
\qed      

	\subsection{Characterization of Global $\varepsilon$-Minimum of Pointwise Supremum}
The intention of this subsection is to introduce the definition of $\varepsilon$-Robust infimum, which appears to be promising  to characterize $\varepsilon$-minimums of the supremum function; we also refer a sufficient condition to guarantee the mentioned property. The introduction of this concept  allows us to cover many ideas reflected in  \emph{max-min theorems} that have been broadly studied, this notation is motivated by the concept of \emph{Robust infimum} or \emph{Decoupled Infimum} used in subdifferential theory to get \emph{fuzzy calculus rules} (see, e.g., \cite{MR2191744,MR2986672,MR3033176,MR2144010}). To be  more precise the reader can observe that when $\varepsilon=0$ the  definition below corresponds to a max-min equation, that is to say, $$\sup_{T}\inf_{X} f_t(x)= \inf_{X}\sup_{T}f_t(x).$$ 
\begin{definition}[$\varepsilon$-Robust infimum]
	We will say that the family of functions $\{f_t : t \in T\}$ has an  $\varepsilon$-Robust infimum on $B\subseteq X$ at $\bar{x}\in B$ provided that 
	\begin{align}\label{epsrobustmin}
	f(\bar{x}) \leq \sup\limits_{t\in T}\inf\limits_{x \in B}f_t(x) +\varepsilon.
	\end{align}
	The special case $\varepsilon=0$ is simply  called a Robust minimum. 
\end{definition}

The min-max problem has been studied in many papers with different states of generality, for this reason it is impossible to recall all the sufficient conditions that establish the interchange  between supremum and infimum; we refer to  \cite{MR1921556,MR2144010,MR791361,MR838482,MR0055678,MR0312194}) for some results and further discussions.  However, we establish the next lemma which guarantees this class of \emph{max-min results} and fits perfectly with the framework of our study.  

First we need the following result.

 {\begin{lemma}\textnormal{\cite[Lemma 4]{MR3507100}} \label{lema6}
	Let  $X$ be a topological space and let $(f_\alpha)_{\alpha \in D}$ be a net  of lsc proper functions defined on $X$ such that
	\begin{equation*}
	\alpha,\beta \in D,\; \alpha \leq \beta \Rightarrow f_\alpha \leq f_\beta.
	\end{equation*}
	For $(\varepsilon_\alpha)_{\alpha \in D} \searrow 0$, let $(x_\alpha)_{\alpha \in D}$ be a relatively compact net such that $$x_\alpha \in \varepsilon_\alpha\text{-}\textnormal{argmin}{f_\alpha} \text{ for each } \alpha. $$  Then
	$$\displaystyle\inf\limits_{x \in X} \sup\limits_{\alpha \in D} f_\alpha(x)= \sup\limits_{\alpha \in D}\inf\limits_{x\in X}f_\alpha(x),
	$$
	and every accumulation point of   $(x_\alpha)$ is a minimizer of the function  $\sup\limits_{\alpha \in D} f_\alpha.$
\end{lemma}
}
\begin{lemma}\label{Lemma:SuficientCondition}\textnormal{[Sufficient condition for robust local minimum]}
	Let $(X,\tau)$ be a topological space and  $B\subseteq X$. Suppose $(T,\preceq)$ is a directed set and the family   of  $\tau$-lsc functions $\{f_t : t \in T\}$ is  increasing, 	$B$ is $\tau$-closed and there exists some $t_0$ such that $f_{t_0}$  is $\tau$-infcompact on B. Then the family $\{f_t : t \in T\}$ has a Robust minimum on $B$.
\end{lemma}
{\it Proof}
	It is easy to see that $\inf\limits_{x\in B}\sup\limits_{t\in T} f_t(x) \geq \sup\limits_{t\in T}\inf\limits_{x \in B}f_t(x)=:\eta$. For the opposite inequality let us assume that $\eta <+\infty$.   Consider  $\varepsilon_t \in (0,1)$ and  points $x_t \in B$ such that $\varepsilon_t\to 0$  and $x_t$ belongs to $\varepsilon_t\text{-}\textnormal{argmin}_{B} \{ f_{t} \}$,  then  $$\inf\limits_{x \in B}f_t(x) \geq \inf\limits_{x \in B}f_{t_0}(x)>-\infty \text{ and } {x_t  \in \{ x :  f_{t_0}(x) \leq \eta + 1  \}} \text{ for all }t \succeq t_0,$$ so the set $(x_t)$ is relatively compact,  and by Lemma \ref{lema6}  we get the result.
\qed   
We finish this section with a very simple result of  the $\varepsilon$-robust infimum in terms of the $\varepsilon$-subdifferential of the initial data.
\begin{proposition}
	Let  $\{f_t : t \in T\}$ be a family of functions which has an  $\varepsilon$-Robust infimum on $B\subseteq X$ at ${x}$, then
	\begin{align}
	0 \in \bigcap\limits_{ \gamma >   0} \bigcup\limits_{t \in T} \partial_{\varepsilon + \gamma} (f_{t}+\delta_B)(x).
	\end{align}
	The condition is also sufficient if $f(x)=f_t(x)$ for all $t\in T$.
\end{proposition}
{\it Proof}
	Let $x$ such that $f(x)\leq \sup_{t\in T}\inf_{y\in B}f_t (y) + \varepsilon$, then for a given $\gamma >0$ one can take $t\in T$ such that $f_t(x)\leq f(x)\leq \inf_{y\in B} f_t(y)+ \varepsilon + \gamma$, that is, $0 \in  \partial_{\varepsilon + \gamma} (f_{t}+\delta_B)(x)$. Conversely, if $f(x)=f_t(x)$ for all $t\in T$ and $0$ belongs to $ \bigcap_{ \gamma >   0} \bigcup_{t \in T} \partial_{\varepsilon + \gamma} (f_{t}+\delta_B)(x)$, one has that for every $\gamma >0$ there exists some $t\in T$ such that  $f(x)=f_t(x) \leq f_t(y) + \varepsilon + \gamma $ for all $y\in B$. Then, $f(x) \leq \sup_{t \in T}\inf_{y\in B} f_t(y) + \varepsilon + \gamma $, so the arbitrariness  of $\gamma$ gives the result.
\qed   

	\section{Pointwise Supremum Function of  an Ordered Family of Epi-pointed Functions}\label{ORDEREDSETS}
	In this section, we investigate the $\varepsilon$-subdifferential of the pointwise supremum under the assumption that the  data is epi-pointed and the set $T$ is a directed set. This extra hypothesis allows us to give better estimations for the conjugate function of the supremum function as well as formulae for the $\limsup$ and $\liminf$ functions. It is worth mentioning that the property of being a directed set can be  satisfied using the family $\mathcal{P}_f(T)$ and the  family of functions  defined by $g_A(\cdot) := \max_{s \in A} f_s$, then $f(x)= \sup_{A \in  \mathcal{P}_f(T)} g_A(\cdot)$. Moreover the epi-pointed property can be also obtained using an appropriate perturbation of the functions $f_t$ as will be used in some of the proofs.
	\begin{theorem}\label{Prop:1}
		Let  $\{f_t : t \in T \}$  be an increasing family  of epi-pointed functions such that $f^\ast$ is proper and  {$f^{\ast \ast} = \sup_{t\in T} f^{\ast\ast}$}.
		Then for every $x^\ast\in \inte \dom f^\ast$ one has
		\begin{align}\label{Equ:Conjungate}
		f^\ast(x^\ast)=\inf\{ f_t^\ast(x^\ast) : t\in T \}.
		\end{align}
	\end{theorem}
	{\it Proof}
	 {Let $x_0^\ast \in\inte \dom f^\ast$,	by Lemma \ref{Lemma2} \ref{Lemma2b} $f^\ast(x_0^\ast)\leq \inf\{ f_t^\ast(x_0^\ast) : t\in T \}$. Now, using Lemma \ref{Lemma2} \ref{Lemma2important2} and \ref{Lemma2d}  we have that $f^\ast$ is the closed hull of the function $\inf\{ f_t^\ast(\cdot) : t\in T \}$, then one can take   $x^\ast_\alpha \to x_0^\ast$  such that \begin{align}\label{limiinf}
		f^\ast(x_0^\ast)= \lim_{\alpha } \inf \{ f_t^\ast(x_\alpha^\ast)  :t \in T	\}, 
		\end{align} besides by Lemma \ref{Lemma2} \ref{Lemma2important}  $x^*_0 \in \inte \dom f_{t_0}$ for some $t_0 \in T$. Because the family of epi-pointed functions $\{ f^\ast_s\}_{s\succeq t_0}$ is decreasing, they are (uniformly) bounded on $x_0^\ast+ \tilde{U}_{t_0}$ (recall that $x^\ast_0 \in \inte \dom f_{t_0}^\ast$ and $f_{t_0}^\ast$ is continuous at $x_0^\ast$), then  by \cite[Theorem 2.2.11]{MR1921556} we can  find a neighbourhood $U_{t_0} \in \mathcal{N}_0$, $K_{t_0} >0$ and a continuous seminorm $\rho$  such that 
		\begin{align}\label{Eq:Lip}
		|f^\ast_s(x^\ast)-f_s^\ast(y^\ast)| \leq K_{t_0} \rho(x^\ast - y^\ast) \text{ for all } x,y \in x^\ast_0  +  U_{t_0}, \text{ for all } s \succeq t_0.
		\end{align}
		Now let $\gamma_\alpha \to 0$ and $t_\alpha \in T$ with $t_\alpha \succeq t_0$  such that $\gamma_\alpha + f^{\ast}(x^\ast_0)  \geq f^{\ast}_{t_\alpha} (x_\alpha^\ast)$ (recall \eqref{limiinf}), then using (\ref{Eq:Lip}) one yields
		$$\gamma_\alpha + f^{\ast}(x^\ast_0)  \geq -K\rho(x_\alpha - x_0^\ast) + f_{t_\alpha}(x_0^\ast) \geq  -K\rho(x_\alpha - x_0^\ast)  + \inf\{ f_t^\ast(x^\ast) : t\in T \}.$$
		Then taking the limits we conclude the result.}
	\qed   
	
	Now based on  Theorem \ref{Prop:1} together with classical calculus rules for the conjugate of the sum of convex functions we give one answer to the question proposed in \cite[Question 3.6]{MR3571567}.
	
	For this purpose we introduce the set $\ell^1(\N,X^\ast)$ as the set of all sequences $(x_n^\ast) \in X^\ast$ such that 
	$\sum_{n\in \N} | \langle  x^\ast_n , x \rangle|  <+\infty $ for every $x\in X$   and the linear functional   $\langle x^\ast , \cdot \rangle  := \sum_{n\in \N} \langle  x^\ast_n , \cdot \rangle $ belongs to $X^\ast$ (see \cite{INTECONV} for more details). 
	
	\begin{corollary}
		Consider  $(x^\ast_n) \in \ell^1(\N,X^\ast)$, $(\alpha_n) \in  \ell^1(\N,\R)$ and  a sequence of  functions $f_n \in \Gamma_0(X)$ such that 
		\begin{align*}
		f_n(x) \geq \langle x^\ast_n , x \rangle  + \alpha_n\;\; \forall x\in X,\; \forall n \in \N.
		\end{align*}
		Consider $f(x):=\sum_{n\in \N} f_n(x)$. If $f$ and $f_n$ are epi-pointed functions, then, for every $x^\ast\in \inte  \dom f^\ast$ one has
		\begin{align}
		f^\ast(x^\ast)&= \inf \{   \sum_{n=1 }^N f_n^\ast(x_n^\ast)  : N\in \N \text{ and } \sum_{n=1 }^N x_n^\ast =x^\ast    \} \label{eqcorollaryinfinitesums01}\\
		&= \inf \{   \sum_{n\in \N } f_n^\ast(x_n^\ast)  :  \sum_{n\in \N} x_n^\ast =x^\ast    \}. \label{eqcorollaryinfinitesums02}
		\end{align}
	\end{corollary}
	{\it Proof}
		We (may) assume $f_n\geq 0$ for all $n\in \N$ (otherwise we can consider $\tilde{f}_n := f_n - x_n^\ast - \alpha_n$). Hence, $f^\ast_n(0) \leq 0$  for all $n\in \N$.	Thus, $f=\sup_{n\in \N} g_n$, where   $g_n:= \sum_{k\leq n} f_k$.  First we recall (see, e.g., \cite[Corollary 2.3.5]{MR1921556}) that for every $ n \in \N$ one has 
		\begin{align}
		\label{eqcorollaryinfinitesums}  (g_n )^\ast(x^\ast)&=  \cl^{w^\ast}\square_{k=1}^n f^\ast_k (x^\ast), \; \forall x^\ast \in X^\ast,\\
		\label{eqcorollaryinfinitesums2} (g_n )^\ast(x^\ast)&=  \square_{k=1}^n f^\ast_k (x^\ast), \; \forall x^\ast \in \inte(\dom   \square_{k=1}^n f^\ast_k),
		\end{align}
		where $\square_{k=1}^n f^\ast_k (\cdot)=  \inf\{ f_1^\ast(x^\ast_1)+...+f_n^\ast(x^\ast_n)  : \sum\limits_{k=1}^n x_k^\ast=\cdot 	\}$.   Now (\ref{eqcorollaryinfinitesums}) implies that $$\sum_{k=1}^n \dom  {f^\ast_k} \subseteq \dom  {g_n^\ast} \subseteq \cl \big( \sum_{k=1}^n \dom  {f^\ast_k} \big),$$
		whence 
		$\inte(\dom   \square_{k=1}^n f^\ast_k)= \inte \dom g^\ast_n$.
		
		Now, consider $x^\ast\in \inte \dom f^\ast$. It is not difficult to prove that  $$f^\ast(x^\ast) \leq \inf \{   \sum_{n\in \N } f_n^\ast(x_n^\ast)  :  \sum_{n\in \N} x_n^\ast =x^\ast    \}.$$  For the other inequality we notice that  $\{ g_n\}$ form an  increasing family of proper lsc convex epi-pointed   functions  and $f=\sup_{n\in \N } g_n $, then using Theorem \ref{Prop:1} one has that
		$f^\ast(x^\ast)=\inf \{ g_n^\ast(x^\ast) : n \in \N  \}$,  besides by  Lemma \ref{Lemma2}  \ref{Lemma2important} $x^\ast\in \bigcup_{n\in \N} \inte \dom g_n^\ast$. Whence   (\ref{eqcorollaryinfinitesums2}) yields 
		$f^\ast(x^\ast)= \inf \{   \sum_{k=1 }^n f_k^\ast(x_k^\ast)  : N\in \N \text{ and } \sum_{k=1 }^n x_k^\ast =x^\ast    \}$, which concludes the proof of (\ref{eqcorollaryinfinitesums01}). Finally, in order to prove (\ref{eqcorollaryinfinitesums02})  we recall $f_n^\ast(0) \leq 0$, thus  (\ref{eqcorollaryinfinitesums01}) implies 
		$f^\ast(x^\ast) \geq  \inf \{   \sum_{n\in \N } f_n^\ast(x_n^\ast)  :  \sum_{n\in \N} x_n^\ast =x^\ast    \}$ and that concludes the proof.
	\qed

				A straightforward application  of Lemma \ref{Lemma2} and Theorem \ref{Prop:1} gives us an estimation for the conjugate of  $\limsup$ and $\liminf$ functions.
				
				\begin{corollary}
					Consider a directed set $T$  {and} an arbitrary family of functions  {$\{f_t : t \in T \} \subseteq  \Gamma_0(X)$}.
					\begin{enumerate}[label={(\alph*)},ref={(\alph*)}]
						\item If $h(x):=\limsup_{t\in T} f_t (x)$  {is proper}, then
						\begin{align}\label{newone}
						\epi h^\ast =\bigcap\limits_{t \in T} \cco\{  \bigcup\limits_{s\succeq t } \epi f_t^\ast   \}.
						\end{align}
						\item If $g(x):=\liminf_{t\in T} f_t(x)$  {is proper}  and there exists $t_0\in T$ such that $(\inf_{s\succeq t_0} f_s))(\cdot)$ is epi-pointed and 	 {$g^{\ast \ast}= \sup\limits_{t\in T} g_t^{\ast \ast},$}
					where 	$g_t:= \inf_{s\succeq t} f_s$, then 
						\begin{align}\label{equation12}
						 {g^\ast}(x^\ast):=\limsup f_t^\ast(x^\ast) \text{ for all } x^\ast  \in \inte \dom  {g^\ast}
						\end{align}
						and  $\epi  {g^\ast} = \bigcup\limits_{t \in T}  		\bigcap\limits_{s\succeq t}	  \epi f_s^\ast$. 
					\end{enumerate}	
				\end{corollary}
				{\it Proof}
					First define $h_t:= \sup_{ s \succeq t} f_s$. Then
					\begin{enumerate}[label={(\alph*)},ref={(\alph*)}]
						\item  {$h^\ast =\sup_{x\in X}\{  \langle \cdot , x \rangle - h(x) \}= \sup_{t\in T}\sup_{x\in X}\{  \langle \cdot , x \rangle - h_t(x)\}=\sup_{t\in T}h_s^\ast$. Then, we notice that  $h_t$ must  be  proper for some $t \in T$, so by  Lemma \ref{Lemma2} \ref{Lemma2important2}   $ \epi h_s^\ast=\cco\{  \bigcup\limits_{s\succeq t } \epi f_t^\ast   \}$, and so \eqref{newone} holds.}
						\item  {$g=\sup_{t\in T} g_t$ and $\{ g_t \}_{t\succeq t_0 }$ is an increasing family of epi-pointed functions  satisfying $g^{\ast \ast}= \sup_{t\in T} g_t^{\ast \ast}$. Then, Theorem  \ref{Prop:1} and Lemma \ref{Lemma2} \ref{Lemma2b} yield (\ref{equation12}), also Lemma \ref{Lemma2} \ref{Lemma2b} implies $\epi g_t^\ast =  		\bigcap_{s\succeq t}	  \epi f_s^\ast$, and consequently $\epi g^\ast =\bigcup_{t \in T}  	\epi g_t^\ast =\bigcup_{t \in T}  		\bigcap_{s\succeq t}	  \epi f_s^\ast$.}
					\end{enumerate}
				\qed
				\begin{theorem}\label{TEO:EpsilonSub}
					Let  $\{f_t : t \in T \}$  be an increasing family of  epi-pointed functions such that  {$f^{\ast \ast}=\sup_{t\in T}f_t^{\ast \ast}$},  then  for all $x\in X$.
					\begin{align}
					\sub_\varepsilon f(x)&=\bigcap\limits_{\substack{t\in T \\ \gamma >0 }} \cl\bigg\{   \bigcup\limits_{ s \succeq   t} \sub_{\varepsilon+\gamma} f_s(x)    \bigg\}\label{FORM:1}\\
				&=	 \bigcap\limits_{\substack{t\in T \\ \gamma >0 }} \cl\bigg\{   \bigcup\limits_{ s_0 \succeq   t} \bigcap_{ s \succeq   s_0}\sub_{\varepsilon+\gamma} f_s(x)    \bigg\} \label{FORM:2}.
					\end{align}
					Moreover, the above formulae hold if the functions $f_t$ are not necessarily \mbox{epi-pointed}, but they belong to $ \Gamma_0(X)$. 
				\end{theorem}
				{\it Proof} First we check that the right side of (\ref{FORM:2}) is included in the $\varepsilon$-subdifferential of $f$ at $x$. Indeed, let $\gamma >0$ and $t\in T$ and pick $x^\ast\in  \bigcup\limits_{ s_0 \succeq   t} \bigcap\limits_{ s \succeq   s_0}\partial_{\varepsilon+\gamma} f_s(x)  $. Then $\langle x^\ast , y - x \rangle \leq f_s(y) -f_s(x) +\varepsilon + \gamma\leq f(y)-f_s(x) +\varepsilon +\gamma$ for all $y \in X$ and  for all $s\succeq s_0$, then $\langle x^\ast , y - x \rangle \leq  f(y)-f(x) +\varepsilon +\gamma$,  {which means $x^\ast \in \sub_{\varepsilon + \gamma} f(x)$}. Therefore 
					$$ \bigcap\limits_{\gamma>0}\bigcap\limits_{\substack{t\in T  }} \cl\bigg\{   \bigcup\limits_{ s_0 \succeq   t} \bigcap_{ s \succeq   s_0}\sub_{\varepsilon+\gamma} f_s(x)    \bigg\} \subseteq \bigcap\limits_{\gamma>0}\sub_{\varepsilon + \gamma} f(x) =\partial_{\varepsilon} f(x).$$
			 {Now, it is easy to see that (\ref{FORM:2}) is included in (\ref{FORM:1}). For the opposite inclusion we notice that for any co-final set $\tilde{T} \subseteq T$ we have
					\begin{align*}
					 \bigcap\limits_{\substack{t\in T \\ \gamma >0 }} \cl\bigg\{   \bigcup\limits_{ s_0 \succeq   t} \bigcap_{ s \succeq   s_0}\sub_{\varepsilon+\gamma} f_s(x)    \bigg\}& = \bigcap\limits_{\substack{t\in \tilde{T} \\ \gamma >0 }} \cl\bigg\{   \bigcup\limits_{ s_0 \succeq   t} \bigcap_{ s \succeq   s_0}\sub_{\varepsilon+\gamma} f_s(x)    \bigg\}.
				\end{align*}
				  Then consider $\gamma >0$ and $t \in T_\gamma(x)$ arbitrary. Pick $s_2 \succeq s_1 \succeq  t$ (it implies that $s_1,s_2 \in T_\gamma(x)$) and  $x^\ast \in  \sub_{\varepsilon+\gamma} f_{s_1}(x)  $, then 
				\begin{align*}
				\langle x^\ast , y - x\rangle&  \leq f_{s_1}(y) - f_{s_1}(x) + \varepsilon+\gamma \\
				&\leq f_{s_2}(y)  - f_{s_1}(x) +f_{s_2}(x)   - f_{s_2}(x)+  \varepsilon+\gamma\\
				&\leq f_{s_2}(y)  - f_{s_1}(x) +f(x)   - f_{s_2}(x) +\varepsilon+\gamma\\
				&\leq f_{s_2}(y)   - f_{s_2}(y)+ \varepsilon+2\gamma,
				\end{align*}
				which means $x^\ast \in  \sub_{\varepsilon+2\gamma} f_{s_2}(x)  $, since it is for every  $s_2 \succeq s_1 \succeq  t$ we conclude that 
				\begin{align*}
				 \bigcup\limits_{ s \succeq   t} \sub_{\varepsilon+\gamma} f_s(x)   \subseteq    \bigcup\limits_{ s_0 \succeq   t} \bigcap_{ s \succeq   s_0}\sub_{\varepsilon+2\gamma} f_s(x),   
				\end{align*}
				and consequently
					\begin{align*}
			\bigcap\limits_{\substack{t\in T \\ \gamma >0 }} \cl\bigg\{   \bigcup\limits_{ s \succeq   t} \sub_{\varepsilon+\gamma} f_s(x)    \bigg\} &\subseteq 	\bigcap\limits_{\substack{t\in T_\gamma(x) \\ \gamma >0 }} \cl\bigg\{   \bigcup\limits_{ s \succeq   t} \sub_{\varepsilon+\gamma} f_s(x)    \bigg\} \\ &\subseteq \bigcap\limits_{\substack{t\in  T_\gamma(x)  \\ \gamma >0 }} \cl\bigg\{   \bigcup\limits_{ s_0 \succeq   t} \bigcap_{ s \succeq   s_0}\sub_{\varepsilon+\gamma} f_s(x)    \bigg\} 
			\\ &\subseteq \bigcap\limits_{\substack{t\in T \\ \gamma >0 }} \cl\bigg\{   \bigcup\limits_{ s_0 \succeq   t} \bigcap_{ s \succeq   s_0}\sub_{\varepsilon+\gamma} f_s(x)    \bigg\} 
				\end{align*}} 
					
					Now we focus on proving that  $\partial_\varepsilon f(x)$ is included in the right side of (\ref{FORM:1}),  {so w.l.o.g. we can assume that $\partial_\varepsilon f(x)\neq \emptyset$, in particular $f^*$ is proper,  and by  Lemma \ref{Lemma2c} $f$ is epi-pointed}. First we prove this in the case that the functions $f_t$ are epi-pointed. Thanks to Lemma \ref{Upper:stimationeta} it is enough to prove that for every $\gamma >0$ and $t\in T$
					
					\begin{align}\label{key1}
					\sub_{\gamma + \varepsilon} f(x) \cap \inte \dom f^\ast \subseteq  \bigcup\limits_{ s \succeq   t}\sub_{\varepsilon+2\gamma} f_s(x)  .
					\end{align}

					Then, take $x^\ast$ in the left side of (\ref{key1}), then  $f^{*}(x^\ast) + f(x) \leq \langle x^\ast ,x \rangle  + \varepsilon+\gamma$, then using Theorem \ref{Prop:1} (recall $x^\ast \in \inte \dom f^\ast$) and due to  the fact that  the family is increasing,  we can take  $s \succeq t$ such that 
					 {$f_s^\ast(x) + f_s(x) \leq \langle x^\ast ,x \rangle  + \varepsilon+2\gamma$}, which implies $x^\ast\in \sub_{\varepsilon+2\gamma} f_s(x)$.
					
					Now when the functions are not necessarily epi-pointed one  can use the following argumentation;  in order to simplify the notation we assume that $x=0$. For the sake of  understanding the proof first consider that $X$ is a \mbox{finite-dimensional} space. Then, for every $\alpha >0$ we define the epi-pointed  function 
					$g^\alpha_t(\cdot) := f_t(\cdot) + \alpha\| \cdot\|$  (see, e.g., \cite[Proposition 3.1]{Correa2017}) and  consider $g^\alpha=\sup_{t\in T} g_t^\alpha = f + \alpha\| \cdot\|$ .  Then, $\sub_\varepsilon f(0) = \bigcap_{\alpha > 0}    \partial_{\varepsilon } g^\alpha (0)$, therefore 
					\begin{align*}
					\partial_\varepsilon  f(0)  &= \bigcap_{\alpha > 0}    \partial_{\varepsilon } g^\alpha (0)=\bigcap_{\alpha > 0}     \bigcap\limits_{\substack{t\in T \\ \gamma >0 }} \cl \bigg\{ \bigcup\limits_{ s \succeq t} \sub_{\varepsilon+\gamma} g_s^\alpha (0)    \bigg\}\\
					&= \bigcap_{\alpha > 0}  \bigcap\limits_{\substack{t\in T \\ \gamma >0 }}   \cl \bigg\{  \bigcup\limits_{ s \succeq t}\big( \sub_{\varepsilon+\gamma} f_s(0) + \alpha \mathbb{B}(0,1)\big) \bigg\}\\	
					&= \bigcap\limits_{\substack{t\in T \\ \gamma >0 }}   \bigcap_{\alpha > 0}  \cl \bigg\{  \big(\bigcup\limits_{ s \succeq t} \sub_{\varepsilon+\gamma} f_s(0)\big) + \alpha \mathbb{B}(0,1) \bigg\}	\\
						&= \bigcap\limits_{\substack{t\in T \\ \gamma >0 }}  \cl \bigg\{  \bigcup\limits_{ s \succeq t} \sub_{\varepsilon+\gamma} f_s(0) \bigg\}
					\end{align*}
						Finally, if the space  $X$ is infinite-dimensional we must identify the elements of 
					$\partial (f +\delta_L)(0)$ with $	\partial f_{|_L} (0)$; here $L$ is a finite-dimensional subspace of $X$ and $f_{|_L}$ is the restriction of $f$ to $L$.
					Indeed, consider  $t\in T$ and $\gamma>0$, then take  $e_k \in X$ for $k=1,...,p$ and  define $V=\{ x^\ast \in X^\ast : |\langle x^\ast, e_k |   \leq 1 \}$ and let $L$ be the linear subspace generated by $\{ e_k \}_{k=1}^p$. Now, we are going to prove that there exists $s \succeq t$ such that $	\partial (f +\delta_L)(0)\subseteq \partial_{\varepsilon+\gamma} f_s(x) + V$. Let $x^\ast \in 	\partial (f +\delta_L)(0)$ and consider $P:X\to L$ a continuous projection and $P^\ast$ its adjoint operator, 
					then $x^\ast_{|_L } \in \partial f_{|_L}(0)$, then by the last part there exists  $s \succeq t$ such that $x^\ast_{|_L } \in \partial (f_s)_{|_L}(0) +W^\ast$,
					where $W^\ast:=\{ x^\ast \in L^\ast : |\langle x^\ast, e_k |   \leq 1 \}$. Then
					$x^\ast = P^\ast(x^\ast_{|_L} ) + x^\ast - P^\ast(x^\ast_{|_L} ) \in  \partial_{\varepsilon+\gamma}( f_s+\delta_{L})(0) +V+  L^\perp$, besides \emph{Hiriart-Urruty and Phelps's formula} (see, e.g., \cite[Theorem 3.2]{MR1330645}) implies $$\partial_{\varepsilon+\gamma}( f_s+\delta_{L})(0) \subseteq \partial_{\varepsilon+\gamma}( f_s+\delta_{L})(0) + V^\ast + L^\perp. $$
					Therefore $x^\ast \in \partial_{\varepsilon+\gamma} f_s(x) + V$.
					
				\qed 
			  
				The following example shows that the closure operator in the above result is necessary even for functions defined in the real line.
				\begin{example}
					Consider the functions $f_n(x)=(1-1/n)| x|$, then $f(x)=|x|$ and $\sub_{\varepsilon} f_n(0)=[-1+1/n,1-1/n]$ and $\sub_{\varepsilon}  f(0)=[-1,1]$ for all $\epsilon \geq 0$.  {Whence $\bigcup_{n\in \N} \sub_{\varepsilon} f_n(0) = (-1,1)\neq [-1,1]=\sub  f(0)$.}
				\end{example}
				The next result corresponds to a general formula for the $\epsilon$-normal set of an intersection of arbitrary closed and convex sets.
			
				\begin{corollary}\label{theoremnormalcone}
					Consider a family of closed convex subsets $\{ C_t : t\in T \}$ and $C:=\bigcap_{t\in T} C_t$. Then for every $\varepsilon \geq 0$ and   {$x\in C$}
					\begin{align}\label{theoremnormalcone:01}
					\textnormal{N}^\varepsilon_C(x) =\bigcap\limits_{ \gamma >   0} \cl \bigg\{ \sum_{ t \in A  }  {\textnormal{N}^{\eta_t}_{C_t}(x)} : \begin{array}{c}
					A \subseteq T,  \#A < +\infty,  \\ \eta_t \geq 0 \text{ and } \sum\limits_{t\in A} \eta_t = \varepsilon+\gamma 
					\end{array} \bigg\}.
					\end{align}
				\end{corollary}
			{\it Proof}
			First the right-hand side of \eqref{theoremnormalcone:01} is included in $\textnormal{N}^\varepsilon_C(x)$. We focus on the opposite inclusion.	Consider the sets $G_A:=\cap_{t\in A} C_t$ for $A \in \mathcal{P}_f(T)$. Then 
				
				$\delta_C = \sup_{A \in\mathcal{P}_f(T)} \delta_{G_A}$, then by Theorem \ref{TEO:EpsilonSub} we have that

				\begin{equation}\label{equation1}
				\textnormal{N}^\varepsilon_C(x)= \bigcap_{ \gamma >   0} \cl\big\{   \textnormal{N}^{\varepsilon+\gamma}_{G_A}(x)  : A \subseteq T, \; \#A <+\infty	\big\},
				\end{equation}
				moreover
				$\textnormal{N}^{\varepsilon+\gamma}_{G_A}(x)= \sub_{\varepsilon+\gamma}\big( \sum_{t\in A}\delta_{C_t}  	 \big)(x)$
				and by Hiriart-Urruty and Phelps's  formula
				\begin{equation}\label{equation2}
				\sub_{\varepsilon+\gamma}\big( \sum_{t\in A}\delta_{C_t}  	 \big)(x)=\cl^{w^\ast} \big\{ \sum_{t\in A}  {\textnormal{N}^{\eta_t}_{C_t}(x) }  : \eta_t \geq 0 \text{ and }\sum_{t\in T} \eta_t = \varepsilon + \gamma    \big\}.  
				\end{equation}
				Consequently, mixing (\ref{equation1}) and (\ref{equation2}) we get the result.
				\qed   
				
			\section{Calculus  for   $\varepsilon$-Subdifferential  and the Fenchel Conjugate  of Pointwise Supremum Function}\label{SECTION:Calculus}
			This section is devoted to the study of the $\varepsilon$-subdifferential  of the pointwise supremum function. In the first part  we  assume that the space $X$ is  \mbox{finite-dimensional}; this condition allows us to use Carathéodory's Theorem to bound the cardinal of the support of the elements $\lambda \in \Delta(T)$ in Proposition \ref{Epsilonformula} and with this we obtain more precise calculus rules for the $\varepsilon$-subdifferential and the Fenchel conjugate of the pointwise supremum function. In the second part of this section,  the  finite-dimensional condition over $X$  is removed using some reduction to  finite-dimensional subspaces; for this reason we prefer to use the notation $X$ in Subsection \ref{FINITEDIMENSION} for a finite-dimensional space also instead of using  the Euclidean space $\R^{n}$, putting  emphasis on the fact that $X$ could be an abstract  finite-dimensional space.
			
			\subsection{Finite-Dimensional Banach Spaces}\label{FINITEDIMENSION}
		In this section we give simplifications of Proposition \ref{Epsilonformula} under some classical \emph{qualification conditions} in a  finite-dimensional Banach space $X$. In this subsection we assume that the functions $f_t$ satisfy the relation  $f^{\ast \ast}=\sup_{t \in T}  f^{\ast \ast}_t$.
			
		Let us introduce the following  qualification condition at a point $x\in X$.
			
			\begin{equation}\label{QC1} 
			\textnormal{N}_{\dom f}(x) \textnormal{ does not contain lines.}
			\end{equation}
	\begin{remark}\label{remarkcondition}
		It is important to recall that for any arbitrary function $f:X\to \Rex$ with proper conjugate the following statements are equivalent $(i)$ $\textnormal{N}_{\dom f}(x)$ does not contains lines, $(ii)$ $\cco f$ is continuous at some point and $(iii)$ $f^\ast$ is epi-pointed. Moreover, $\big[ \sub_{\varepsilon} f(x) \big]_{\infty} =\textnormal{N}_{\dom f}(x)$ for all $\varepsilon \geq 0$ and all $x\in X$ such that $\sub_\varepsilon f(x) \neq \emptyset$.
		\end{remark}
	
	The key tool to establish this exact formulation is the \emph{pointedness} of the normal cones involved and  Carathedory's Theorem, which allow us  to give a limiting representation of (\ref{FORMULABASIS}) and (\ref{FORMULABASIS2}). The next lemma represents the major ideas in the proofs of the main results for this section.

	\begin{lemma}\label{Lema51}
		Consider $\varepsilon \geq 0$, $n:=\min \{ \#T, \dim X + 1\}$ and   $x^\ast\in \partial_\varepsilon f(x)$, then  there are sequences $t_{i,k} \in T$, $(\lambda_{i,k}) \in \Delta(\{1,...,n\})$,  $x^\ast_{i,k} \in \partial_{\eta_{i,k} } f_{t_{i,k}} (x)$ and $\eta_{i,k}\geq 0$  with  $i\in \{1,...,n \}$, $k\in \N$ such that $\sum_{i=1}^n \lim\limits_{k}  (\lambda_{i,k} \cdot \eta_{i,k})  \leq \varepsilon$,  
		\begin{align*}
		\sum\limits_{i=1}^{n} \lim_{k} \lambda_{i,k} \cdot ( f(x)-f_{t_{i,k}}(x))  \leq \varepsilon  - \sum\limits_{i=1}^{n} \lim_{k} \lambda_{i,k} \eta_{i,k}
		\end{align*}
		
		and $\sum_{i=1}^n \lambda_{i,k} \cdot x_{i,k}^\ast \to x^\ast$. Moreover, one of the following conditions holds.
		\begin{enumerate}[label={(\alph*)},ref={(\alph*)}]
			\item\label{parta} There exists $n_1 \in \N$ with $n_1 \leq  n $ such that $\lambda_{i,k} \overset{k \to \infty}{\longrightarrow} \lambda_i > 0$, $x_{i,k} \overset{k \to \infty}{\longrightarrow} x_i^\ast$, $\eta_{i,k} \to \eta_i$ for $i \leq n_1$   and $\lambda_{i,k} \overset{k \to \infty}{\longrightarrow}  0$, $\lambda_{i,k}\cdot x_{i,k} \overset{k \to \infty}{\longrightarrow} x_i^\ast$ for $n_1 < i \leq n$, 
			\begin{align*}
			\sum\limits_{i=1}^{n_1} \lambda_i (\lim\limits_{k \to \infty}( f_{t_{i,k}}(x) -f(x))  + \sum\limits_{i>n_1} \lim\limits_{k \to \infty}(\lambda_{i,k}  f_{t_{i,k}}(x) -f(x)) \\\leq \varepsilon  - \sum\limits_{i=1}^{n_1} \lambda_i \lim\limits_{k \to \infty} \eta_{i,k} -\sum\limits_{i>n_1} \lim\limits_{k \to \infty}  \lambda_{i,k}\eta_{i,k},
			\end{align*}
	$	\sum\limits_{i=1}^{n_1} \lambda_i =1$	and  $x^\ast=\sum\limits_{i=1}^{n_1} \lambda_i x^\ast_i + \sum\limits_{i>n_1}^{n} x_i^\ast$, or
			\item \label{Lema51Partb} There are $\nu_k \searrow 0$ such that $\nu_k \cdot \lambda_{i,k}\cdot x^\ast_{i,k} \overset{k \to \infty}{\longrightarrow} x_i^\ast$ and $\sum\limits_{i=1}^{n_1}x^\ast_i =0$ with not all $x_i^\ast$ equal to zero.
			
			Moreover, if one assumes  (\ref{QC1}), then  \ref{parta} always holds.
		\end{enumerate}
	\end{lemma}
	{\it Proof} Consider $x^\ast \in \partial_{\varepsilon} f(x)$ and $\gamma_k \to 0 $. Then by  Proposition \ref{Epsilonformula} there are sequences $t_{i,k} \in T$, $(\lambda_{i,k}) \in \Delta(\{1,...,n\})$, $\eta_{i,k}\geq 0$ and $x_{i,k} \in \partial_{\eta_{i,k} } f_{t_{i,k}} (x)$   such that $\sum_{t\in T}\lambda_{t_{i,k}} \cdot \eta_{t_{i,k}} \in [0,\varepsilon +\gamma_k]$
					and 	$\sum_{t\in T} \lambda_t \cdot  f_t  (x) \geq f(x) + \sum_{t\in T}\lambda_t \cdot \eta_t - \gamma $ and $u_k^\ast:=\sum\limits_{i=1}^{n} \lambda_{i,k} \cdot x_{i,k}^\ast \to x^\ast$. We may assume (up to a subsequence)  that  for every $i=1,...,n$  $(i)$ $\lambda_{i,k} \to \lambda_i$, $(ii)$ $\lim\limits_{k}  (\lambda_{i,k} \cdot \eta_{i,k})$  exists and $(iii)$ $\lim_{k} \lambda_{i,k} \cdot ( f(x)-f_{t_{i,k}}(x)) $ exists. 
				Then,  $\sum_{i=1}^n \lim\limits_{k}  (\lambda_{i,k} \cdot \eta_{i,k})  \leq \varepsilon$ and 
				\begin{align*}
				\sum\limits_{i=1}^{n} \lim_{k} \lambda_{i,k} \cdot ( f(x)-f_{t_{i,k}}(x))  \leq \varepsilon  - \sum\limits_{i=1}^{n} \lim_{k} \lambda_{i,k} \eta_{i,k}.
				\end{align*}
					  Moreover (reordering if it is necessary),
				we  assume that $\lambda_{i} >0$ for  $i=1,...,n_1$ and $\lambda_{i} >0$ for $i > n_1$.
				
					Now, on the one hand, if we assume that $\sup_{k, i} \| \lambda_{i,k} \cdot x^\ast_{i,k}\|_{*}  <+\infty$, then up to a subsequence one can assume  that $\lambda_{i,k} \overset{k \to \infty}{\longrightarrow} \lambda_i > 0$, $x_{i,k} \overset{k \to \infty}{\longrightarrow} x_i^\ast$, for $i \leq n_1$   and $\lambda_{i,k} \overset{k \to \infty}{\longrightarrow}  0$, $\lambda_{i,k}\cdot x_{i,k} \overset{k \to \infty}{\longrightarrow} x_i^\ast$ for $n_1 < i \leq n$, 
					\begin{align*}
					\sum\limits_{i=1}^{n_1} \lambda_i (\lim\limits_{k \to \infty}( f_{t_{i,k}}(x) -f(x))  + \sum\limits_{i>n_1} \lim\limits_{k \to \infty}(\lambda_{i,k}  f_{t_{i,k}}(x) -f(x)) \\
					\leq \varepsilon  - \sum\limits_{i=1}^{n_1} \lambda_i \lim\limits_{k \to \infty} \eta_{i,k} + \sum\limits_{i>n_1} \lim\limits_{k \to \infty}  \lambda_{i,k}\eta_{i,k}
					\end{align*}
					and  $x^\ast=\sum\limits_{i=1}^{n_1} \lambda_i x^\ast_i + \sum\limits_{i>n_1}^{n} x_i^\ast$. 
					
					On the other hand, if $\sup_{k, i} \| \lambda_{i,k} \cdot x^\ast_{i,k}\|_{*} =\sup\limits_{k }\{\max\limits_{i} \| \lambda_{i,k} \cdot x^\ast_{i,k}\|    \}=+\infty$ up to a  subsequence one can assume that $\nu_k:=(\max\limits_{i} \| \lambda_{i,k} \cdot x^\ast_{i,k}\|_{*} )^{-1} \searrow 0 $ and $\nu_k \lambda_{i,k}\cdot x^\ast_{i,k} \overset{k \to \infty}{\longrightarrow} x_i^\ast$  with not all $x_i^\ast$ equal to zero, thus $\nu_k u_k^\ast \to  \sum_{i=1}^n x_i^\ast= 0$.
					
					Finally, it is not difficult to see that \ref{Lema51Partb} contradicts (\ref{QC1}).

				\qed   
			Now we present the main result of this section.
			\begin{theorem}\label{TEO5.2}
				Assume that (\ref{QC1}) holds at $x\in \dom f$. Then for every  $\varepsilon \geq 0$ one has
				\begin{align*}
				\partial_{\varepsilon}f(x) &=\bigcup \bigg\{ S(x,\varepsilon_1)+\textnormal{N}^{\varepsilon_2}_{\dom f}(x) \;\big|\;	(\varepsilon_1,\varepsilon_2) \in \Delta^\varepsilon(\{1,2\})	\bigg\},\\
				&= \bigcup \bigg\{ 	\overline{S}(x,\varepsilon_1)+\textnormal{N}^{\varepsilon_2}_{\dom f}(x) \;\big|\; (\varepsilon_1,\varepsilon_2) \in \Delta^\varepsilon(\{1,2\})	\bigg\}, 
				\end{align*}
				where
				\begin{align*}
				S(x,\varepsilon_1)  &:=\bigcap\limits_{\gamma >0 } \cl \bigg\{  \sum_{t\in \supp \lambda } \lambda_t  \partial_{\frac{\varepsilon_t}{\lambda_t} + \gamma } f_t(x)  : \begin{array}{c}
				\lambda \in \Delta(T), \; 
				(\varepsilon_t)\in \Delta^{\varepsilon_1} (T) \text{ and}\\
			f_t(x) +\varepsilon_t/\lambda_t +\gamma \geq 	f(x)     
				\end{array}   \bigg\}\\
				\overline{S}(x,\varepsilon_1)  &:=\bigg\{  \sum_{t\in \supp \lambda } \lambda_t x_t^\ast  : \begin{array}{c}
				\lambda \in \Delta(T), 
				 \;  (\varepsilon_t) \in \Delta^{\varepsilon_1} (T) \text{ and}\\
					x^\ast_t \in A(x,\varepsilon_t,\lambda_t)   
				\end{array}   \bigg\}\\
				\end{align*} 
				and $ A(x,\varepsilon,\lambda):=\bigg\{ x^\ast \in X^\ast :	\begin{array}{c}
				 \exists t_k \in T, \; x_{k}^\ast \in \partial_{\frac{\varepsilon}{\lambda} + \gamma_k}  f_{t_k}(x), \; \gamma_k \searrow 	0,\\ \text{such that } x_k^\ast \to x^\ast \text{ and } \lim f_{t_k}(x) + \varepsilon/\lambda  \geq f(x) 
								\end{array}  \bigg\}$

			\end{theorem}
{\it Proof}
	We focus on the nontrivial  inclusions.  Let $x^\ast \in \partial_\varepsilon f(x)$, then consider sequences   $t_{i,k} \in T$, $(\lambda_{i,k}) \in \Delta(\{1,...,n\})$, $\eta_{i,k}\geq 0$ and $x_{i,k} \in \partial_{\eta_{i,k} } f_{t_{i,k}} (x)$    as in  Lemma \ref{Lema51}, and by (\ref{QC1})  condition  \ref{parta}  must hold. Now, using the notation of Lemma \ref{Lema51} \ref{parta}, define 
	$y^\ast:= \sum_{i=1}^{n_1} \lambda_{i} x^\ast_i $, $z^\ast:= \sum_{i=n_1+1}^n x^\ast_i$, $\varepsilon_1:= 	\sum_{i=1}^{n_1} p_i$, with $p_i:= \lambda_i ( f(x)-\lim\limits_{k \to \infty}f_{t_{i,k}}(x)   +  \eta_{i})$ for  $i=1,..., n_1$ and 
	\begin{align}\label{key11}
	\varepsilon_2:=	  \sum\limits_{i>n_1} \lim\limits_{k \to \infty}\lambda_{i,k} ( f(x) -f_{t_{i,k}}(x))  + \sum\limits_{i>n_1} \lim\limits_{k \to \infty}  \lambda_{i,k}\eta_{i,k},
	\end{align}
	so $x^\ast=y^\ast+z^\ast$ and $\varepsilon_1+\varepsilon_2\leq \varepsilon$. It follows that $y^\ast \in \overline{S}(x,\varepsilon_1)$. Indeed,  define  {$\gamma_{i,k}= |\frac{p_i}{\lambda_i} - (f(x)-f_{t_{i,k}}(x)   +  \eta_{i,k})|$}, then one gets $x^\ast_{i,k} \in \partial_{\frac{p_i}{\lambda_i} + \gamma_{i,k}} f_{t_{i,k}}(x)$ and $\lim_{k} f_{t_{i,k}}(x) +p_i/\lambda_i = f(x) + \eta_i \geq f(x)$, which means  $x^\ast_i \in A(x,p_i,\lambda_i)$.
	
	  Now we show that $z^\ast \in  \textnormal{N}^{\varepsilon_2}_{\dom f}(x)$ and  $y^\ast \in 	{S}(x,\varepsilon_1)$. Let $\gamma$ any number in  $(0, \min\{ \lambda_i : i=1,...,n_1 \})$ and define  $M= \sup\{ \|  x^\ast_{i,k}\|_{*} : k\in \N, \; i=1,...,n_1 \}$, let $k \in\N$ be such that $\| y^\ast - \sum_{i=1}^{n_1}  \lambda_{i,k} x^\ast_{i,k}\|_{*} \leq \gamma$, $(1- \sum_{i=1}^{n_1} \lambda_{i,k}) M < \gamma$, 
	\begin{align}\label{key7.1}
	  \lambda_{i,k} (f(x)- f_{t_{i,k}}(x)  +  \eta_{i,k} )< p_i +\gamma^2
	\end{align}

	and $\varepsilon_1 ( \frac{1 }{ \lambda_{i,k}}- \frac{1}{ (1- \sum_{i=1}^{n_1} \lambda_{i,k})  + \lambda_{i,k} } )< \gamma$ for all $i =1,...,n_1$. Then, we set $\lambda \in \Delta(T)$, $(\varepsilon_t) \in \Delta^{\varepsilon_1}(T)$ and $\eta \in \R_{+}^{(T)}$ by 
	
	\begin{align*}
	\lambda_t =& \left\{ 
	\begin{array}{ccl}
	\lambda_{i,k} + (1- \sum_{i=1}^{n_1} \lambda_{i,k}), & & \text{if } t=t_{1,k},\\
	\lambda_{i,k}, &  & \text{if } t=t_{i,k} \text{ for }i=2,...,n_1,\\
	0, & &\text{ otherwise},
	\end{array}\right.\\
	\varepsilon_t =& \left\{ 
	\begin{array}{ccl}
	p_i, & & \text{if } t=t_{i,k},\\
	0, & &\text{ otherwise},
	\end{array}\right.\\
	\eta_t =& \left\{ 
	\begin{array}{ccl}
	\eta_{i,k},  & & \text{if } t=t_{i,k} \text{ for }i=1,...,n_1,\\
	0, & &\text{ otherwise},
	\end{array}\right.
	\end{align*}
	respectively. Then $y^\ast \in  \sum\limits_{t\in \supp \lambda } \lambda_t \partial_{\eta_t} f_t(x) + \mathbb{B}(0,2\gamma)$, moreover  (recall (\ref{key7.1}) and $\gamma < \lambda_{i,k}$) for $t\neq t_{1,k}$,  $f_t(x) +\varepsilon_t/\lambda_t +\gamma \geq 	f(x)    +\eta_t $, and for $t= t_{1,k}$, $f_t(x) +\varepsilon_t/\lambda_{1,k} +\gamma \geq 	f(x)    +\eta_t $, besides
	  $$ \frac{\varepsilon_t}{\lambda_{1,k}} = \frac{\varepsilon_t}{\lambda_t }  + \varepsilon_t \big(    \frac{1}{\lambda_{1,k} } - \frac{1}{\lambda_{1,k} + (1- \sum_{i=1}^{n_1} \lambda_{i,k})} \big) \leq  \frac{\varepsilon_t}{\lambda_t }  + \gamma.$$
	 Therefore $y^\ast \in  \sum\limits_{t\in \supp \lambda } \lambda_t \partial_{\frac{\varepsilon_t}{\lambda_t} + 2\gamma } f_t(x) + \mathbb{B}(0,2\gamma)$
	
	Finally, for every $y\in \dom f$, one has
	\begin{align*}
	\langle z^\ast , y- x\rangle &= \lim\limits_{k \to \infty} \sum_{i=n_1+1}^n \lambda_{i,k} \langle  x^\ast_{i,k} ,  y-x \rangle \\
	&\leq \lim\limits_{k \to \infty} \sum_{i=n_1+1}^n \lambda_{i,k} \left( 			f_{t_{i,k}} (y) - f_{t_{i,k}}(x)	 + \eta_{i,k}\right) \\
	&\leq \lim\limits_{k \to \infty} \sum_{i=n_1+1}^n \lambda_{i,k} ( f(y) - f(x) )  +	\lim\limits_{k \to \infty} \sum_{i=n_1+1}^n \lambda_{i,k} ( f(x)   - f_{t_{i,k}} (x))\\
	&\leq \varepsilon_2 \;\;(\textnormal{recall } (\ref{key11})).
	\end{align*}
\qed   
Although the assumption that the normal does not contain  lines, or some of its equivalent statements (see Remark \ref{remarkcondition}), is standard in convex analysis, for some applications it is better to understand this condition in terms of the functions $f_t$, for this reason we introduce the next  \emph{qualification condition} in terms of the normal cones of the data function $\{ f_t\}_{t\in T}$. 
\begin{equation}\label{QC2}
\begin{array}{c}
\text{For every } A \subseteq \mathcal{P}_f(T)  \text{ one has }\\
x^\ast_t \in \textnormal{N}_{\dom f_t } (x) \text{ for } t\in A \text{ and }  \sum_{t \in A} x^\ast_t = 0 \Longrightarrow  x_t^\ast =0 \text{ for all } t\in A.
\end{array}
\end{equation}
Under the additional assumptions of the compactness of $T$ and some \emph{continuity property} of the function $t \to f_t(w)$ we can prove that (\ref{QC2}) is equivalent to (\ref{QC1}). More precisely we get the following result.
\begin{theorem}\label{teorem:ref}
	 Assume that $T$ is a compact space, $t \to f_t(w)$ is upper semicontinous for every $w \in X$ and \eqref{QC2} holds at $x$. 
	Then, for every $\varepsilon \geq 0$
	\begin{align}
	\textnormal{N}^\varepsilon_{\epi f}(x,f(x))&=  \bigcup\bigg\{   \sum_{t\in\supp (\varepsilon_t)} \textnormal{N}^{\varepsilon_t}_{\epi f_t}(x,f(x)): 
	\begin{array}{c}
(\varepsilon_t) \in \Delta^{\varepsilon}(T)
	\end{array}	 \bigg\}\label{Lemmanormalcone1}\\
	\partial_{\varepsilon }f(x) &=\bigcup \bigg\{	\mathcal{S}(x,\varepsilon_1, T) + \mathcal{N}(x,\varepsilon_2,T) : (\varepsilon_1,\varepsilon_2) \in \Delta^\varepsilon(\{1,2\}) \bigg\} \label{teorem:ref:eq1} 
	\end{align}
	where 
	\begin{align*}
	\mathcal{S}(x,\varepsilon_1)& := \bigg\{   \sum_{t\in\supp \lambda} \lambda_t \partial_{\varepsilon_t/\lambda_t} f_t(x) : \begin{array}{c}
	\lambda \in \Delta(T),\; (\varepsilon_t) \in \Delta^{\varepsilon_1}(T) \\
\text{ and } f_t(x) +\varepsilon_t/\lambda_t  \geq 	f(x)  \end{array}  \bigg\},\\
	\mathcal{N}(x,\varepsilon_2)&:=  \bigg\{   \sum_{t\in\supp \eta} \textnormal{N}^{\eta_t}_{\dom f_t} (x) : 
	(\eta_t) \in \Delta^{\varepsilon_2}(T) \bigg\}.
	\end{align*}
\end{theorem}
{\it Proof}   First we prove (\ref{Lemmanormalcone1}), the result is trivial if $f^{\ast \ast}=-\infty$, or $f^{\ast \ast}=+\infty$. Therefore, we assume that $f^{\ast \ast}$ is proper. Consider the set $\tilde{T}:=\{ t : f^{\ast \ast}_t \neq-\infty\}$
and the functions $g_t := \delta_{\epi f_t}$ and $g :=\sup_{t\in \tilde{T}} g_t$.  The relation $f^{\ast \ast}=\sup f^{\ast \ast}$ is equivalent to $g^{\ast \ast}=\sup_{t\in \tilde{T}} g^\ast_t$. Then, we pick  $(x^\ast,\alpha)$ in $\textnormal{N}^\varepsilon_{\epi f}(x,f(x))=\partial_\varepsilon g (x,f(x))$, so  applying Lemma \ref{Lema51} (and following its notation) there are sequences $t_{i,k} \in T$, $(\lambda_{i,k}) \in \Delta(\{1,...,n\})$, $\eta_{i,k}\geq 0$ and elements $(x^\ast_{i,k},\alpha_{i,k} )$ in $\partial_{\eta_{i,k}} g_{t_{i,k}} (x,f(x))=\textnormal{N}^{\eta_{i,k}}_{\epi  f_{t_{i,k}} } (x,f(x))$   with  $i\in \{1,...,n \}$, $k\in \N$ such that $\sum_{i=1}^n \lim\limits_{k}  (\lambda_{i,k} \cdot \eta_{i,k})  \leq \varepsilon$, 
and $\sum_{i=1}^n \lambda_{i,k} \cdot (x_{i,k}^\ast,\alpha_{i,k}) \to (x^\ast,\alpha)$. We recall that necessarily $\alpha \leq 0$ and  $\alpha_{i,k} \leq 0$ for all $i,k$.  Since $T$ is compact we may assume (up to a subnet) that $t_{i,k} \to t_i$, let us define $\varepsilon_t=\sum\limits_{i: \; t_i =t}\lim\limits_{k}  (\lambda_{i,k} \cdot \eta_{i,k})$ if there exists some $t_i=t$, and $\varepsilon_t=0$ if $t\neq t_i$ for all $i$.

  Now suppose  that condition \ref{Lema51Partb} of Lemma \ref{Lema51}  holds, then the elements $\nu_k \lambda_{i,k} \cdot (x_{i,k}^\ast,\alpha_{i,k})$ converges to  $ (x^\ast_i,\beta_i)$ and $\sum_{i=1}^n (x^\ast_i,\beta_i) =(0,0)$ for some $\nu_k \searrow 0$ and not all  $ (x^\ast_i,\beta_i) $ are equal to zero, then necessarily $\beta_i =0$, because $\beta_i \leq 0$ for all $i=1,...,n$. Now we check that  $x_i^\ast \in  \textnormal{N}_{\dom f_{t_i}}(x)$, indeed let $y \in \dom f_{t_i}$, because $t \to f_t(y)$ is upper semicontinuous (usc) we have that there exists $r \in \R$ such that $(y,r) \in \epi{f_{t_{i,k}}}$ for large enough $k$, consequently
  \begin{align*}
  \langle x^\ast_i, y - x  \rangle &= \lim_{k}\big(  \langle   \nu_k    \lambda_{i,k} \cdot x_{i,k}^\ast, y -x \rangle +  \nu_k   \lambda_{i,k}  \alpha_{i,k} (r - f(x)) \big)\\
  &= \lim_{k} \nu_k    \lambda_{i,k} \big(  \langle   x_{i,k}^\ast, y -x \rangle +  \alpha_{i,k} (r - f(x)) \big)\\
  &\leq  \lim_{k}    \nu_k \lambda_{i,k} \eta_{i,k} \leq  \lim_{k}    \nu_k \varepsilon =0,
  \end{align*}
the last means $x_i^\ast \in  \textnormal{N}_{\dom f_{t_i}}(x)$ and $\sum_{i=1}^n x_i^\ast=0$ with not all $x^\ast_i$ equal to zero,  this contradicts (\ref{QC2}). Therefore, condition \ref{parta} of Lemma \ref{Lema51} must hold,  then $\lambda_{i,k} \cdot (x_{i,k}^\ast,\alpha_{i,k})  \to (x^\ast_i,\beta_i)$ and  using that  $t \to f_t(w)$ is usc we get 
  \begin{align*}
  (u^\ast_t,\alpha_t ):=\sum_{j:\ t_j=t} (x^\ast_i,\beta_i) & \in  \textnormal{N}^{\varepsilon_t}_{\epi{f_t}} (x,f(x)),  \text{ if there exists some  } t_i=t,\\
  (u^\ast_t,\alpha_t ):=(0,0)& \in \textnormal{N}^{\varepsilon_t}_{\epi{f_t}} (x,f(x)),  \text{ otherwise}, \\
  \end{align*}
  and  $\sum_{t \in T} (u^\ast_t,\alpha_t )= (x^\ast,\alpha)$. Consequently using (\ref{Lemmanormalcone1}) (with $\varepsilon=0$) we get that (\ref{QC2}) implies  (\ref{QC1}). Then, we can  apply Theorem \ref{TEO5.2} and following its notation we use the compactness of $T$ and the upper semicontinuity of $t \to f_t(w)$ (together with similar arguments as in the first part)  to prove  that $\overline{S}(x,\varepsilon_1) $ is contained in $\mathcal{S}(x,\varepsilon_1)$. Moreover, using the fact that $(x^\ast,0 )\in \textnormal{N}^{\varepsilon_2}_{\epi f}(x,f(x))$  iff $x^\ast\in \textnormal{N}^{\varepsilon_2}_{\dom f}(x)$   together with  (\ref{Lemmanormalcone1}) we derive that $\textnormal{N}^{\varepsilon_2}_{\dom f} (x)$ is equal to $ \mathcal{N}(x,\varepsilon_2)$. This concludes the proof of (\ref{teorem:ref:eq1}). 
\qed   
\begin{remark}
	It has not escaped our notice that using (\ref{teorem:ref:eq1}) one can prove that 
	\begin{align}
	\partial_{\varepsilon} f(x) &=	\bigcup \bigg\{ 		\mathcal{S}(x,\varepsilon_1, T_1) + \mathcal{N}(x,\varepsilon_2,T_2) : \begin{array}{c}
	(\varepsilon_1,\varepsilon_2) \in \Delta^\varepsilon(\{1,2\}),	\\T_1 \cap T_2 = \emptyset \text{ and}  \\
	\#T_1+ \# T_2 \leq \textnormal{dim}(X)+1
	\end{array} 						\bigg\},\label{teorem:ref:eq2}
	\end{align}
where 
\begin{align*}
	\mathcal{S}(x,\varepsilon_1, T_1)& := \bigg\{   \sum_{t\in\supp \lambda} \lambda_t \partial_{\varepsilon_t/\lambda_t} f_t(x) : \begin{array}{c}
		\lambda \in \Delta(T),\; (\varepsilon_t) \in \Delta^{\varepsilon_1}(T_1) \\
		\text{ and } f_t(x) +\varepsilon_t/\lambda_t  \geq 	f(x)  \end{array}  \bigg\},\\
	\mathcal{N}(x,\varepsilon_2,T_2) &:=  \bigg\{   \sum_{t\in\supp \eta} \textnormal{N}^{\eta_t}_{\dom f_t} (x) : 
	(\eta_t) \in \Delta^{\varepsilon_2}(T_2) \bigg\}.
\end{align*}
	Indeed, to prove  (\ref{teorem:ref:eq2}) we notice that by the finite dimension of $X$ every element  $x^\ast \in \sub_\varepsilon f(x)$ must be expressed as $x^\ast= \sum_{t\in T_1 } \lambda_ t x^\ast_t  + \sum_{s \in T_2} w_s^\ast$ with some  $x^\ast_t \in   \partial_{\varepsilon_t/\lambda_t} f_t(x)$, $t\in T_1$, $w^\ast_s \in  \textnormal{N}^{\nu_s}_{\dom f_s} $, with $\sum_{t\in T_1} \varepsilon_t + \sum_{s\in T_2} \nu_s \leq \varepsilon$,  $T_1,T_2 \subseteq T$,  $\#(T_1 \cup T_2) \leq \dim X+1$ and $f_t(x) +\varepsilon_t/\lambda_t  \geq 	f(x) $ for all $t \in T_1$. Then $\lambda_t \partial_{\varepsilon_t/\lambda_t} f_t(x) +  \textnormal{N}^{\nu_t}_{\dom f_t}(x) \subseteq \lambda_t \partial_{(\varepsilon_t+ \nu_t)/\lambda_t } f_t(x) $ for $t \in T_1\cap T_2$, consequently $x^\ast \in 	\mathcal{N}(x,\varepsilon_1 ,T_1) + 	\mathcal{N}(x, \varepsilon_2,T_2\backslash T_1)$, where $\varepsilon_1:= \sum_{t\in T_1} \varepsilon_t  + \sum_{t\in T_1\cap T_2}\nu_t $ and $\varepsilon_2:=   \sum_{s\in T_2\backslash T_1} \nu_s $.
	\end{remark}

The final goal of this section is to give the following characterization of the epigraph of  $f^\ast$ and consequently an expression for itself. This result can be understood as the conjugate counterpart of Theorems \ref{TEO5.2} and \ref{teorem:ref}.
\begin{theorem}\label{THEOREMCONJUGATE}
	
	\begin{enumerate}[label={(\alph*)},ref={(\alph*)}]
			\item If $f^\ast$ is epi-pointed, one has 
		\begin{align*}
		\epi f^\ast &= \co \cl \bigcup\{ \epi f_t: t \in T \} + (\epi f^\ast)_\infty.\\
		\end{align*}
	 Consequently $f^\ast = \co \ell  +(f^\ast)^{\infty}$, where $\ell$ is the lsc function such that its epigraph is  $\cl \bigcup\{ \epi f_t: t \in T \}$.
		\item\label{THEOREMCONJUGATEcaseb} If $T$ is compact and $t \to f_t(x)$ is upper semicontinous for every $x\in X$ and for every $B \in \mathcal{P}_f(T)$  and every $(x^\ast_t,\alpha_t) \in (\epi f_t^\ast)_{\infty}$ with $t\in B$ one has 
		\begin{align}\label{Infiniteconstraintcompactcase2}
	\sum_{t \in B} (x^\ast_t,\alpha_t) = (0,0) \Longrightarrow  (x_t^\ast,\alpha_t) =(0,0) \text{ for all } t\in B.
		\end{align}	
	Then
		\begin{align*}
		\epi f^\ast &=\co \{  \epi f^\ast_t  : t\in T\} + \co\{  (\epi f_t^\ast)_\infty : t\in T \}.
		\end{align*}
	\end{enumerate}
\end{theorem} 
{\it Proof}
	Assume that $f^\ast$ is epi-pointed. Consider $(x^\ast,\alpha) \in \epi f^\ast$, then by Lemma \ref{Lemma2} \ref{Lemma2important2}  we have $(x^\ast,\alpha)=\lim (x^\ast_n,\alpha_n)$ for some $(x^\ast_n,\alpha_n) \in \co\{  \bigcup_{t\in T}\epi f^\ast_t \} $. Moreover,  we can write 
	$$(x^\ast_n,\alpha_n) = \sum_{k=1}^N \lambda_{t_{k,n}} (y^\ast_{k,n},\beta_{k,n})$$ where $N:=\dim X+2$, $(\lambda_{t_{k,n}} ) \in \Delta(T)$ and $(y^\ast_{k,n},\beta_{k,n}) \in \epi f^\ast_{t_{k,n}}$ for some ${t_{k,n}} \in T$. Now by a similar argumentation as the one given in Lemma \ref{Lema51}, up to a subsequence  $\lambda_{t_{k,n}} (y^\ast_{k,n},\beta_{k,n})$ must converge to a point $(w^\ast_k,\beta_k)$ and $\lambda_{k,n} \to \lambda_{k}$, otherwise there exists $(u^\ast_k,\mu_k) \in (\epi f^\ast)_{\infty}=\epi (f^\ast)^\infty$ not all equal to zero such that $\sum_{k=1}^{N} (u_k^\ast,\mu_k)=(0,0)$ which contradicts the fact that $f^*$ is epi-pointed.  Now, on the one hand for every $k =1,...,N$ with $\lambda_k\neq 0$ one has that $$(y^\ast_{k,n},\beta_{k,n}) \to (y^\ast_k , \beta_k) \in \cl \bigcup\{ \epi f_t: t \in T \} .$$  On the other hand for every $k =1,...,N$ with $\lambda_k= 0$  one has that $$\lambda_{k,n} (y^\ast_{k,n},\beta_{k,n})  \to (w^\ast_k, \gamma_k) \in (\epi f^\ast)_{\infty}.$$ Therefore, $$x^\ast = \sum\limits_{k: \lambda_k>0}\lambda_k (y^\ast_k , \beta_k) + (w^\ast,\beta),$$ where $ (y^\ast_k , \beta_k)  \in \cl \bigcup\{ \epi f_t: t \in T \}$ and  $(w^\ast,\beta)=\sum\limits_{k: \lambda_k=0} (w^\ast_k, \gamma_k) \ \in(\epi f^\ast)_\infty$. 
	The proof of   \ref{THEOREMCONJUGATEcaseb} follows
	 similar arguments as the previous one, together with the argumentation about the compactness of $T$ and the upper semicontinuity as was given in the proof of Theorem \ref{teorem:ref}, so we omit the proof.
\qed   
\begin{remark}
	A relevant point to consider is that all the convex combinations given in the results of this subsection can be  written as  convex combinations of no more than $\dim X +1 $ points (or $\dim X+2$  for formulas related to epigraphs). In the same form,  condition (\ref{QC2}) and (\ref{Infiniteconstraintcompactcase2}) are  equivalent to  consider $A \subseteq T$ with $\#A \leq \dim X$ and $B\subseteq T$ with  $\#B \leq \dim X +1$ in  (\ref{QC2}) and (\ref{Infiniteconstraintcompactcase2}) respectively.
	\end{remark}
\subsection{Infinite-Dimensional Locally Convex Spaces}\label{INFINITEDIMENSION}

For the sake of simplicity  we assume that the functions $f_t\in \Gamma_0(X)$ for all $t\in T$. 

Consider $x\in X$,  by $\mathcal{F}_x$ we denote the family of all finite-dimensional linear subspaces of $X$  containing  $x$.

\begin{theorem}\label{TEOREMINFDIMENSION}
	Consider a  family of  functions $\{ f_t\}_{t\in T}\subseteq \Gamma_0(X)$. Then, for every $x\in X$ and $\varepsilon \geq 0$ one has
	\begin{align}
	\partial_{\varepsilon}f(x)&=\bigcap\limits_{\substack{L\in  \mathcal{F}_x\\ \gamma >0 }}\cl^{w^\ast}\bigcup \bigg\{ 	\mathfrak{S}(x,\varepsilon_1,\gamma)+\textnormal{N}^{\varepsilon_2}_{\dom f\cap L}(x) \;\big|\;	(\varepsilon_1, \varepsilon_2)\in \Delta^\varepsilon(\{1,2\})	\bigg\},\label{EQ3:TEOREMINFDIMENSION}
	\end{align}
	where
	\begin{align*}
	\mathfrak{S}(x,\varepsilon_1,\gamma)&:= \bigg\{  \sum_{t\in T} \lambda_t  \partial_{\frac{\varepsilon_t}{\lambda_t}+\gamma} f_t(x)  : \begin{array}{c}
	\lambda \in \Delta(T), \;	(\varepsilon_t) \in \Delta^{\varepsilon_1}(T), \\
		\text{ and }
	f_t(x) +\varepsilon_t/\lambda_t +\gamma \geq 	f(x)  
	\end{array}   \bigg\}
	\end{align*} 
\end{theorem}
{\it Proof} 
	Since the right side of the equation is contained in the $\varepsilon$-subdifferential of $f$ at $x$,  we must focus on the opposite inclusion. Without loss of generality we (may) assume that $x=0$. 
	
	Now take  $\gamma >0$, 	$V^\ast \in  \mathcal{N}_0$  and $L\in \mathcal{F}_0$ such that $L^{\perp} \subseteq V^\ast$. Take $x^\ast$ in $\partial_{\varepsilon} f(0)$, consider $F:= \textnormal{spn}\{ L\cap \dom f\} $	and   $P:X \to F$ a continuous linear projection and let us denote by $P^\ast$ its adjoint operator. Then $x^\ast = y^\ast + z^\ast$, where $y^\ast:=P^\ast(x_{|_F}^\ast)$ and $z^\ast: =x^\ast - y^\ast\in F^{\perp}$,  then $x_{|_F}^\ast $ belongs to $ \partial f_{|_F}(0)$. Since the relative interior of $\dom f_{|_F} $ with respect to $F$ is nonempty, the hypotheses of  Theorem \ref{TEO5.2} hold (in $F$), then there exists $\varepsilon_1,\varepsilon_2 \geq 0$ with $\varepsilon_1+\varepsilon_2=\varepsilon$ together with elements $u^\ast$,  $\lambda^\ast$ in $S(0,\varepsilon_1)$ and  $ \textnormal{N}^{\varepsilon_2}_{\dom f_{|_F}} (0)$ respectively such that $x_{|_F}^\ast=u^\ast + \lambda^\ast$, which implies the existence of $	\lambda \in \Delta(T)$, $(\varepsilon_t) \in \Delta^{\varepsilon_1}(T)$ such that $$f_t(x) +\varepsilon_t/\lambda_t +\gamma \geq 	f(x) \text{ and  }u^\ast \in \sum_{t \in \supp \lambda}\partial_{\frac{\varepsilon_t}{\lambda_t}+\gamma} (f_t)_{|_F}(x) +(P^\ast)^{-1}(V^\ast).$$ Whence $P^\ast(u^\ast) \in \sum_{t \in \supp \lambda}\partial_{\frac{\varepsilon_t}{\lambda_t}+\gamma} (f_t+ \delta_{F})(x) + V^\ast$ and by \cite[Theorem 3.2]{MR1330645} we have 
	\begin{align*}
	\sum_{t \in \supp \lambda}\partial_{\frac{\varepsilon_t}{\lambda_t}+\gamma} (f_t+ \delta_{F})(x)&=\cl^{w^\ast}\big( 	\sum_{t \in \supp \lambda}\partial_{\frac{\varepsilon_t}{\lambda_t}+\gamma} f_t(0) + L^\perp \big)\\
	&\subseteq \partial_{\frac{\varepsilon_t}{\lambda_t}+\gamma} f_t(0) +V^\ast.
	\end{align*}

	Moreover,   $P^\ast(\lambda^\ast) + z^\ast \in N_{\dom f \cap L}(x)$. Therefore, 
	$$x^\ast \in 	\mathfrak{S}(x,\varepsilon_1,\gamma) + \textnormal{N}^{\varepsilon_2}_{\dom f\cap L}(x) +V^\ast+V^\ast.$$
	From the arbitrariness  of  $L$, $V^\ast$ and $\gamma>0$ we get the result.
\qed   
\begin{remark}
	In \cite{MR2960099} the author has proposed changing the family $\mathcal{F}_x$ for some family  of convex sets $C \subseteq \dom f$ which covers $\dom f$. It is not difficult to see that $\textnormal{N}^\varepsilon_{\dom f  \cap L}(x) \subseteq \textnormal{N}^\varepsilon_{\dom f \cap C \cap L}(x)$ for every  set $C$ containing $x$. Then, in order to get  	(\ref{EQ3:TEOREMINFDIMENSION}) with some special class of sets $C$ one needs a formula in the following form
	\begin{equation}\label{normal1}
	\textnormal{N}^\eta_{\dom f \cap C \cap L}(x) = \cl^{w^\ast} \big(	\textnormal{N}^{\eta}_{\dom f \cap C }(x) +(L-x)^\perp	 \big).
	\end{equation}
	or some $\nu$-enlargement of the above expression 
		\begin{equation}\label{normal2}
	\textnormal{N}^\eta_{\dom f \cap C \cap L}(x) = \bigcap\limits_{ \nu >   0}\cl^{w^\ast} \big(	\textnormal{N}^{\eta+\nu}_{\dom f \cap C }(x) +(L-x)^\perp	 \big).
	\end{equation}
	Then, in order to obtain a new formula for the $\varepsilon$-subdifferential of  $f$ one can consider a family of sets $\{ C_\alpha : \alpha \in \mathcal{A} \}$ satisfying  (\ref{normal1}) for large enough $L$  and all $\eta \leq \varepsilon$ which cover $\dom f$, or if one prefers  the use of the  $\nu$-enlargement, one can consider a family   of sets $\{ C_\alpha : \alpha \in \mathcal{A} \}$ satisfying  (\ref{normal2}) for large enough $L$ and all $\eta \leq \varepsilon$ which cover $\dom f$. In particular it happens in \cite[Theorem 6.2]{MR2960099}, where the author has used a family of closed convex sets . 
	 
	 A more concrete idea is to understand the $\varepsilon$-normal set of the intersection as the $\varepsilon$-subdifferential of the sum of the indicators of the respective sets, that is,
	\begin{equation}\label{key30}
	\textnormal{N}^{\eta}_{\dom f \cap C \cap L} (x) =\left\{  \begin{array}{l}
	\partial_{\eta} (\delta_{ \dom f \cap C \cap L })(x) = \partial_{\eta} (\delta_{ \dom f\cap C} +\delta_L)(x),  \text{ if } \eta >0,\\
	\textnormal{N}_{\R_{+}(\dom f)  \cap L} (x)= \partial (\delta_{ \R_{+} (\dom f\cap C -x) +x} +\delta_L )(0),   \text{ if } \eta =0,
	\end{array}\right.
	\end{equation}
then one can  apply  Hiriart-Urruty and Phelps's formulae. Nevertheless, these formulae need  that the sets be closed, or at least satisfied
\begin{align}
\cl(\dom f\cap C \cap L)&=  \cl(\dom f\cap C) \cap L\label{normal3}\\
\cl(  (\R_{+} (\dom f\cap C -x) +x) \cap L)&= \cl({ \R_{+} (\dom f\cap C -x) +x}) \cap L\label{normal4}
\end{align}
  (see, e.g., \cite{INTECONV},\cite[Proposition 2]{MR3522006}, where   Hiriart-Urruty and Phelps's formulae have been imposed with weaker assumptions, in particular (\ref{normal1})  and (\ref{normal2}) under  (\ref{normal3})  and (\ref{normal4})).  Then, in order to obtain formulae for the $\varepsilon$-subdifferential (with $\varepsilon>0$) of $f$ one defines  $\mathcal{A}_x^\varepsilon$ as some family of sets $C$ satisfying (\ref{normal3}) and (\ref{normal4})  for large enough $L$. Similarly, for the \emph{Moreau-Rockafellar} Subdifferential one can consider $\mathcal{A}_x$ as some family of sets $C$ satisfying  (\ref{normal4})  for large enough $L$. 
  
  Typical examples of sets that satisfy (\ref{normal3}) and (\ref{normal4}) are sets which satisfy the \emph{accessibility lemma}, that is there exist some $x_0 \in \dom f\cap C$ such that  $$[x_0, x[ :=\{ (1-\lambda) x_0 + \lambda x: \lambda \in [0,1]	\} \subseteq \dom f\cap C \text{ for every } x \in \cl(\dom f\cap C).$$ This assumption is satisfied under the nonemptiness of \emph{some notations of relative interior} (see \cite{MR1992991}).  In this scenario, we can underline \cite[Corollary 8]{MR2448918} where the authors assume that the \emph{relative interior} of the domain is not empty, in this case it is enough to use $\mathcal{A}_x =\{ \dom f  \}$. Another assumption is that   $\R_{+} (\dom -x)$ is closed (see, e.g., \cite[Corollary 9]{MR3561780}), again it is enough to consider $\mathcal{A}_x =\{ \dom f  \}$.

\end{remark}
As a  corollary to our formula we recover  the result given in  \cite[Theorem 4]{MR2448918} (see also \cite{MR2489616,MR3561780,MR2837551}).
\begin{corollary}\label{COROLLARYCONTI}
	In the setting of  Theorem \ref{TEOREMINFDIMENSION} one has for all $x\in X$
	\begin{align}
	\sub f(x) &=\bigcap\limits_{\substack{L\in  \mathcal{F}_x\\ \gamma >0 }}\cl^{w^\ast}\bigg\{ \co\big\{ \bigcup\limits_{t \in T_\gamma(x) }  \sub_{\gamma } f_t(x)\big\}   +\textnormal{N}_{\dom f\cap L}(x) \bigg\},\label{EQ3.1:TEOREMINFDIMENSION}
	\end{align}
\end{corollary}

\begin{corollary}
	In the setting of Theorem \ref{TEOREMINFDIMENSION} assume that $f$ is continuous at some point of its domain, then for $\varepsilon\geq 0$ and  all $x \in X$ one has
	\begin{align*}
	\sub_{\varepsilon}f(x) &=\bigcup \bigg\{ \bigcap\limits_{\gamma >0 } \cl\mathfrak{S}(x,\varepsilon_1+\gamma,\gamma) +\textnormal{N}^{\varepsilon_2}_{\dom f}(x) \;\big|\;	\varepsilon_1 +\varepsilon_2=\varepsilon \text{ and } \varepsilon_1,\varepsilon_2\geq 0	\bigg\},
	\end{align*}
\end{corollary}
{\it Proof}
Let $x^\ast \in \sub_\varepsilon f(x)$ and $\gamma >0$, then by  Theorem \ref{TEOREMINFDIMENSION} there exist  nets $y^\ast_{\nu,L} \in \mathfrak{S}(x,\varepsilon_1(\nu,L),\gamma)$ and $\lambda^\ast_{\nu,L} \in \textnormal{N}^{\varepsilon_2(\nu,L)}_{\dom \cap L}(x)$ such that $x^\ast= \lim y^\ast_{\nu,L}  + \lambda^\ast_{\nu,L}$; because $f$ is continuous at some point the net $(y^\ast_{\nu,L} )$ must be bounded, so (up to a subnet) we may assume $y^\ast_{\nu,L}  \to y^\ast$ and $\lambda^\ast_{\nu,L} \to \lambda^\ast$ and $\varepsilon_1(\nu,L) \to \varepsilon_1$ and $\varepsilon_2(\nu,L) \to \varepsilon_2$. It follows that  $\lambda^\ast \in \textnormal{N}^{\varepsilon_2}_{\dom f}(x)$. Now it only remains to show that $y^\ast \in \cl\mathfrak{S}(x,\varepsilon_1+\gamma,\gamma)$. Indeed, consider $(\nu_0,L_0)$ such that  for every $(\nu,L) \succeq (\nu_0,L)$, $|\varepsilon_1(\nu,L) - \varepsilon_1| \leq \gamma$. Whence $y^\ast_{\nu,L} \in \sum_{t\in \supp \lambda} \lambda_t \sub_{\varepsilon_t/\lambda_t + \gamma } f_t (x)$ and  $f_t(x) + \varepsilon_t/\lambda_t + \gamma \geq f(x)$ for some $(\varepsilon_t) \in \Delta^{\varepsilon(\nu,L)}(T)$, then considering the real numbers $\tilde{\varepsilon}_t:= \varepsilon_t + \lambda_t (\gamma + \varepsilon_1-  \varepsilon_1(\nu,L) $ we get  $y^\ast_{\nu,L} \in \mathfrak{S}(x,\varepsilon_1+\gamma,\gamma)$ for every $(\nu,L) \succeq (\nu_0,L)$, that concludes the proof.
\qed

Due to the fact that our formulae involve the $\varepsilon$-normal set of $\dom f \cap L$ at a point $x$, we present the following result, where the set is expressed  in terms of the data function $f_t$. The proof of the following lemma  follows  \cite[Lemma 5]{MR2448918}.
\begin{lemma}\label{lemma:normalcone}
	In the setting of Theorem \ref{TEOREMINFDIMENSION} assume that $f$ is proper. Then for every $x\in \dom I_f$ and every $\varepsilon \geq 0$, we have that:
	\begin{align}
	\textnormal{N}^{\varepsilon}_{\dom f}(x)
	&=\{ x^\ast \in X^\ast :   (x^\ast,\langle x^\ast,x\rangle+\varepsilon) \in  \epi (\sigma_{\dom f})    \}\label{lemma:normalconeeq:1}\\
	&=\{ x^\ast \in X^\ast :   (x^\ast,\langle x^\ast,x\rangle+\varepsilon) \in  (\epi f^\ast)_{\infty}    \}\label{lemma:normalconeeq:2}\\
	&=\{ x^\ast \in X^\ast :   (x^\ast,\langle x^\ast,x\rangle+\varepsilon) \in  \bigg[ \cco\big( \bigcup_{ t\in T} \epi f^\ast_t:\big)\biggr]_{\infty} \}\label{lemma:normalconeeq:3}\\
	&=\{ x^\ast \in X^\ast :  (x^\ast,\langle x^\ast,x\rangle+\varepsilon) \in \begin{array}{c}
 \bigg[ \cco\big( \bigcup_{ t\in T} \grafo f^\ast_t  \big) \biggr]_{\infty} \\+ \{0 \}\times[0,\varepsilon] \label{lemma:normalconeeq:4}
	\end{array}\}.
	\end{align} 
	Consequently for every $L\in \mathcal{F}_x$,  $	\textnormal{N}^{\varepsilon}_{\dom f \cap L}(x)$ can be expressed as
		\begin{align}
	\{ x^\ast \in X^\ast :   (x^\ast,\langle x^\ast,x\rangle+\varepsilon) \in  \bigg[ \cco\big(  L^\perp\times \R_{+} \cup  \bigcup_{t\in T} \epi f^\ast_t \big)\biggr]_{\infty} \}, \label{lemma:normalconeeq:5}\text{ or}\\
	\{ x^\ast \in X^\ast :  (x^\ast,\langle x^\ast,x\rangle+\varepsilon) \in \begin{array}{c}
	\bigg[ \cco\big(  L^\perp\times \{0\} \cup \bigcup_{t \in T} \grafo f^\ast_t \big) \biggr]_{\infty} \\+ \{0 \}\times[0,\varepsilon] 
	\end{array}\}.\label{lemma:normalconeeq:6}
	\end{align} 
	 
\end{lemma}
{\it Proof}
	 The  equalities \eqref{lemma:normalconeeq:1} and  \eqref{lemma:normalconeeq:2}  follow from the definition of $\textnormal{N}^{\varepsilon}_{\dom I_f}(x)$ and from the fact that $(\epi f^\ast)_{\infty} =(\epi \sigma_{\dom f}) $ (see \cite[Lemma 5]{MR2448918}). The third equation is given by the fact that  $\epi f^\ast= \cco \{ f_t^\ast : t\in T\}$ (see Lemma   \ref{Lemma2}). Now we  have to prove \eqref{lemma:normalconeeq:4}, from the fact that $$\cco\big( \bigcup_{t \in T} \epi f_t^\ast \big) \supseteq \cco \big( \bigcup_{t\in T} \grafo f^\ast_t\big) + \{0\}\times [0,\varepsilon]$$  the inclusion $\supseteq$  holds in \eqref{lemma:normalconeeq:4}, then we   only have to prove the opposite inclusion. We have that 
	\begin{align*}
	\cl^{w^\ast} \bigg[  \cco \big( \bigcup_{t\in T}   \grafo f^\ast_t \big)  +\{0 \}\times\R_{+}\bigg]= \cco \big( \bigcup_{t\in T}    \epi f_t^\ast \big) .
	\end{align*}
	Since $f^\ast$ is proper, we have 
$$	\begin{array}{l}
	\big[ \cco \big( \bigcup_{t\in T}   \grafo f^\ast_t \big)    \big]_{\infty}\cap - \big[\{0 \}\times\R_{+}\big]_\infty \\
	\\
	\subseteq \big[ \cco \big( \bigcup_{t\in T}    \epi f_t^\ast \big) \big]_{\infty}\cap - \big[\{0 \}\times\R_{+}\big]_\infty=\{(0,0)\},
	\end{array} $$
	so  Dieudonné's Theorem (see \cite[Theorem 2.3.1]{MR562254}) implies that $$\cco \big( \bigcup\{   \grafo f^\ast_t: t\in T \} \big)  +\{0 \}\times\R_{+}$$ is closed, besides $ \cco\big(  \bigcup_{ t\in T} \epi f_t^\ast\big)= \cco \big(  \bigcup_{ t\in T} \grafo f^\ast_t \big)  +\{0 \}\times\R_{+}$ and 
	\begin{align}\label{equation:asymptotic}
	\big[ \cco\big(  \bigcup_{ t\in T} \epi f_t^\ast\big) \big]_\infty = \big[\cco \big(  \bigcup_{ t\in T} \grafo f^\ast_t \big) \big]_\infty + \{0 \}\times\R_{+}.
	\end{align}
	Then take $x^\ast \in X^\ast$ such that   $(x^\ast,\langle x^\ast,x\rangle+\varepsilon) \in  \big[ \cco\big(  \bigcup_{ t\in T} \epi f_t^\ast\big) \big]_\infty $, by (\ref{equation:asymptotic}) there exist $(y^\ast,\gamma)\in \big[\cco \big(  \bigcup_{ t\in T} \grafo f^\ast_t \big) \big]_\infty$ and $\eta \geq 0$ such that $(x^\ast,\langle x^\ast,x\rangle+\varepsilon)=(y^\ast,\gamma+\eta)$, so $x^\ast=y^\ast$. Moreover, since 
	\begin{align*}
	\dom f\times\{-1\}\subseteq \bigg[ \big[ \epi f^\ast\big]_{\infty}\bigg]^{-}&=\bigg[\big[ \cco\big(  \bigcup_{ t\in T} \epi f_t^\ast\big) \big]_\infty \bigg]^{-}\\
&	\subseteq \bigg[\big[\cco \big(  \bigcup_{ t\in T} \grafo f^\ast_t \big) \big]_\infty \bigg]^{-},
	\end{align*}
 we get  $\langle (x^\ast, \gamma), (x,-1)\rangle \leq 0$, so $\eta \leq \varepsilon$, and then $$(x^\ast,\langle x^\ast,x\rangle+\varepsilon) \in \big[\cco\{ \grafo f^\ast_t : t\in T\} ]_{\infty} +  \{0 \}\times[0,\varepsilon].$$
 
 Finally, (\ref{lemma:normalconeeq:5}) and (\ref{lemma:normalconeeq:6}) follow applying (\ref{lemma:normalconeeq:3}) and (\ref{lemma:normalconeeq:4})  to the family of functions $\{f_t : t \in T\} \cup\{ \delta_L\}$.
 
\qed   
\section{Conclusions}
We have derived formulae for the $\varepsilon$-subdifferential and the conjugate functions of the supremum function. We studied two different methods to get these expressions. 

The first one corresponds to the case when the index set is a directed set and the family of functions  is an increasing family of epi-pointed functions. Nevertheless, these  hypotheses can be obtained, as we have shown in Theorem   \ref{TEO:EpsilonSub} and Corollary \ref{theoremnormalcone}, considering the maximum function over all the finite parts of the index set and doing some appropriate perturbation of the function in order to get the epi-pointedness property, in this case the formula obtained does not involve the use of normal cones of the supremum function, which for example, allows us to give a formula for the normal cone of an arbitrary intersection only considering the normal cones of the data, which (at least not directly) cannot be derived from the formulae expressed in Section \ref{SECTION:Calculus}. 

The second one corresponds to obtaining general formulae for the conjugate  and the $\epsilon$-subdifferential  of $f$ in \mbox{finite-dimensional} spaces. Next, we have shown how to derive these results in general locally convex spaces by a reduction to  finite-dimensional subspaces. It is important to say that this reduction to  finite-dimensional subspaces allows us to use standard arguments in subdifferential theory to get   the results, which in particular extend some of the most general formulae to the \emph{Moreau-Rockafellar} subdifferential of the supremum function, to the calculus of the $\varepsilon$-subdifferential of the supremum function.

\bibliographystyle{spmpsci_unsrt}
\bibliography{bibliography}

\end{document}